\documentclass[12pt, twoside, reqno]{amsart}
\usepackage[left=0.9in,right=0.9in,top=1.2in,bottom=1.2in]{geometry}
\usepackage[utf8]{inputenc}
\usepackage[T1]{fontenc}
\usepackage{amsmath,amssymb,latexsym,amsthm,mathrsfs}
\usepackage{caption}
\usepackage{tikz-cd}
\usepackage{mathtools}
\usepackage[pdftex]{hyperref}
\usepackage[noabbrev,capitalize]{cleveref}
\usepackage{aliascnt}
\usepackage{tikz}
\usepackage[dvipsnames]{xcolor}

\theoremstyle{plain}
\newtheorem{theorem}{Theorem}[section]
\newaliascnt{lemma}{theorem}
\newtheorem{lemma}[lemma]{Lemma}
\aliascntresetthe{lemma}
\newaliascnt{proposition}{theorem}
\newtheorem{proposition}[proposition]{Proposition}
\aliascntresetthe{proposition}
\newaliascnt{corollary}{theorem}
\newtheorem{corollary}[corollary]{Corollary}
\aliascntresetthe{corollary}
\newaliascnt{definition}{theorem}
\newtheorem{definition}[definition]{Definition}
\aliascntresetthe{definition}

\newtheoremstyle{claim}%
  {3pt}%
  {3pt}%
  {}%
  {0pt}%
  {\itshape}%
  {.}%
  {.5em}%
  {\thmname{#1}\thmnumber{ #2}\thmnote{ (#3)}}%
\theoremstyle{claim}
\newaliascnt{remark}{theorem}
\newtheorem{remark}[remark]{Remark}
\aliascntresetthe{remark}
\newaliascnt{example}{theorem}
\newtheorem{example}[example]{Example}
\aliascntresetthe{example}

\theoremstyle{plain}
\newtheorem{TheoremA}{Theorem}

\newcommand{\proofend}{\hspace*{\fill}~$\Box$\\}

\def\s{\smallskip}
\def\im{\operatorname {im}}
\def\Spin{\operatorname{Spin}}
\def\Int{\operatorname{int}}
\def\Ham{\operatorname{Ham}}
\def\GL{\operatorname{GL}}
\def\Fix{\operatorname{Fix}}
\def\Hom{\operatorname{Hom}}

\def\cl{{\mathcal L}}
\def\cp{{\mathcal P}}
\def\ct{{\mathcal T}}

\def\C{\mathbb{C}}
\def\N{\mathbb{N}}
\def\Q{\mathbb{Q}}
\def\R{\mathbb{R}}
\def\T{\mathbb{T}}
\def\Z{\mathbb{Z}}

\def\CP{\C P}

\def\SO{\operatorname{SO}}

\def\pp{\partial}
\def\inv{^{-1}}
\def\Ch{\rm{Ch}}
\def\sth{\,\vert\,}

\def\pkm{\cp_{k,m}}

\def\ni{\noindent}

\newcommand{\sslash}{\mathord{/\mkern-6mu/}}

\def\proof{\noindent {\it Proof. \;}}
\newcommand{\proofof}[1]{\ni {\it Proof of #1. }}

\usetikzlibrary{calc}
\usepackage{enumitem}

\crefname{TheoremA}{Theorem}{Theorems}
\Crefname{TheoremA}{Theorem}{Theorems}
\usetikzlibrary{arrows.meta}
\newcommand{\cT}{\mathcal{T}}
\definecolor{AppOuter}{RGB}{23,100,123}
\definecolor{AppInner}{RGB}{135,83,155}
\definecolor{AppZero}{RGB}{57,129,94}
\definecolor{AppQuarter}{RGB}{196,115,34}

\begin{document}

\title[Lagrangian tori in the three-dimensional quadric and mirror symmetry]{Lagrangian tori in the three-dimensional quadric and mirror symmetry}

\author{Jo\'e Brendel}  

\address{
D-MATH,
ETH Zürich, 
Rämistrasse 101,
8092 Zürich,
Switzerland }
\email{joe.brendel@math.ethz.ch}

\author{Andrea Piccirilli}

\address{
\parbox[t]{\dimexpr\linewidth-\parindent\relax}{
Institut de math\'ematiques,
Universit\'e de Neuchâtel,
Rue Emile-Argand 11,
\newline
\qquad 2000 \mbox{Neuchâtel}, Switzerland
\vspace{0.2 cm}}}
\email{andrea.piccirilli@unine.ch}

\date{\today}
\maketitle

\begin{center}
    \emph{With an appendix joint with} Soham Chanda
\end{center}

\begin{abstract}
In the three-dimensional quadric, we study the Oakley--Usher torus and its Chekanov counterpart as defined by the authors in the companion paper \cite{OU1}. Computing the versal deformation of the $\Psi$-invariant of Shelukhin--Tonkonog--Vianna, we find that the Chekanov Oakley--Usher torus is knotted in a strong sense: it is not Hamiltonian isotopic to the Biran-lift of any Vianna torus in the two-dimensional quadric. Furthermore, we investigate enumerative properties of the Oakley--Usher torus in $Q^3$ and prove its non-displaceability, settling an open question asked by Oakley--Usher. We relate this to mirror symmetry by proving a split-generation result for the monotone Fukaya category of $Q^3$. This should be viewed as the analogue of a result by Abouzaid--Diogo after compactifying. Interestingly, in the compactification, the monotone copy of the $(S^1 \times S^2)$-generators of Abouzaid--Diogo becomes displaceable. On the other hand, the Oakley--Usher torus still appears as a generator.

In the appendix joint with Soham Chanda, we give an alternative construction of the Chekanov Oakley--Usher torus using a double solid torus mutation. This allows us to compute its disk-potential using wall-crossing results. 
\end{abstract}

\section{Introduction}

In the companion paper~\cite{OU1}, the authors have constructed knotted versions of the Oakley--Usher Lagrangians in \cite{OakUsh16}. The goal of this paper is to have an in-depth look at the situation in the three-dimensional quadric $Q^3$.

\subsection{Setup}
We equip~$Q^3=\{z_0^2+\cdots+z_4^2=0\}\subset\CP^4$ with the Fubini--Study form normalized so that a line has area~$2\pi$. By \cite[Proposition~3.10]{OU1}, this space can be viewed as the symplectic cut of a the cotangent bundle $T^*S^3$. This identifies~$D^*S^3$ with the complement of a copy $Q^2$ of the two-dimensional quadric. We write~$\Sigma\subset Q^3$ for the image of the zero-section under the symplectic cut; in this model~$\Sigma=\{[1:iq]:q\in S^3\}$.

The families~$\cp_{1,1}(r)$ and~$\cp_{1,1}^{\Ch}(a)$ arise by lifting curves under the~$T^2 = \SO(2)\times\SO(2)$ symplectic reduction described in Section~3 in~\cite{OU1}. The ordinary Oakley--Usher Lagrangians $\cp_{1,1}(r)$ are lifts of the noncontractible circles, while their Chekanov counterparts~$\cp_{1,1}^{\Ch}(a)$ are lifts of contractible circles bounding area~$a/2$ in the smooth reduced cylinder. We use the size convention of~\cite{OU1}: both families have their monotone member at~$r=a=4\pi/3$. Throughout this paper we abbreviate
\[
\cp=\cp_{1,1}(4\pi/3),\qquad
\cp^{\Ch}=\cp_{1,1}^{\Ch}(4\pi/3).
\]
Both are Lagrangian three-tori. They are Lagrangian isotopic by~\cite[Theorem B]{OU1}, but not Hamiltonian isotopic, as we shall see now.

\subsection{Main results} As discussed in \cite[Theorem 1.3]{OakUsh16}, the torus $\cp$ can be viewed as the Biran circle bundle lift of the Clifford torus in~$Q^2\cong S^2\times S^2$~\cite[Theorem~1.3]{OakUsh16}. Similarly, one can lift all Vianna tori from~$Q^2$~\cite{Via17}. Our first result shows that the Chekanov OU-torus is different from tori obtained in this way.

\begin{TheoremA}
	\label{thm:quadric_exotic}
	The torus $\cp^{\rm Ch}$ cannot be mapped by a symplectomorphism to $\cp$ or the lift of any other Vianna torus from the divisor $Q^2 \subset Q^3$.
\end{TheoremA}

The proof computes the versal deformations of the~$\Psi$-invariant of Shelukhin--Tonkonog--Vianna~\cite{SheTonVia24}. To that end, we study a Gelfand--Cetlin fibration which has $\cp^{\Ch}$ as its monotone fibre. Related lifting constructions are considered in earlier work by Chanda--Hirschi--Wang \cite{ChaHirWan24} and by Brendel \cite{Bre23b}. This approach will be significantly generalized in upcoming work by Diogo--Tonkonog--Vianna--Wu~\cite{DTVW26}.

Our second result concerns $\cp$. It computes the signed Maslov two disk-potential and identifies a finite set of generators for the monotone Fukaya category.

\begin{TheoremA}
	\label{thm:quadric_non_displ}
    There is a spin structure on the Oakley--Usher torus $\cp\subset Q^3$, and a basis of $\pi_1(\cp) \cong \Z^3$ for which the superpotential $\mathscr{W}_{\mathcal{P}} : \Hom(\pi_1(\mathcal{P}),\Lambda^{\times})\cong (\Lambda^{\times})^3 \to \Lambda$ is
    \begin{equation*}
    \mathscr{W}_{\mathcal{P}}(x_1,x_2,x_3)
    =
    T^{2\pi/3} \left(x_1+x_2+\frac{x_3}{x_1}+\frac{x_3}{x_2}+\frac{1}{x_3}
\right),
    \end{equation*}
    where~$\Lambda$ is the universal Novikov field over~$\C$ defined in~\eqref{eq:NovikovField}. Furthermore, it carries three critical local systems, which, when taken together with a disjoint Lagrangian sphere $\Sigma \subset Q^3$, split-generate the monotone Fukaya category of $Q^3$.
\end{TheoremA}

In particular, the Oakley--Usher torus $\cp$ is non-displaceable, answering a question raised in~\cite{OakUsh16} for this torus. The sphere~$\Sigma$ accounts for the zero-eigenvalue summand of quantum multiplication by~$c_1(Q^3)$, while the three critical local systems on~$\cp$ account for the three nonzero eigensummands. The proof of \cref{thm:quadric_non_displ} relies on a toric degeneration of~$Q^3$ and on the disk-counting results of Nishinou--Nohara--Ueda~\cite{NisNohUed10,NisNohUed12}.

\subsection{Comparison to previous work}
The split-generation result in \cref{thm:quadric_non_displ} is inspired by work of Abouzaid--Diogo \cite{AboDio23}, which treats the open case of $T^*S^3$. Let us briefly compare the two results. Abouzaid--Diogo prove that the monotone Fukaya category is split-generated by the zero section equipped with an uncoutable family of bounding cochains, together with either of the following two families: 
\begin{enumerate}
    \item a one-parameter family of copies of Lagrangian $S^1 \times S^2$ corresponding to the Oakley--Usher Lagrangians $\cp_{0,2}(r) \subset T^*S^3$, equipped with an uncountbale family of suitable local systems; 
    \item a one-parameter family of Oakley--Usher tori $\cp_{1,1}(r) \subset T^*S^3$, equipped with an uncountbale family of suitable local systems.
\end{enumerate}
See \cite{OakUsh16} or \cite{OU1} for more details on these families. After compactification, only one member of each of the families above is monotone in $Q^3$. Interestingly, the monotone copy of $S^1 \times S^2$ is displaceable in $Q^3$ and therefore vanishes in the Fukaya category. The monotone copy of the tori is $\cp \subset Q^3$ and thus appears as a generator by our result.

Kim~\cite{Kim23} and Kawamoto~\cite{Kaw24} study a different torus in $Q^3$. In fact, their torus can be viewed as the lift of the Chekanov torus in~$Q^2$. Enumerative properties of the OU-torus in~$T^*S^3$ were studied by Chan--Pomerleano--Ueda~\cite{CPU16} and Achig-Andrango~\cite{Ach24}; the latter also discusses the effect of compactification on the potential. Sheridan~\cite{She16} established homological mirror symmetry for Fano hypersurfaces (including $Q^3$) in complex projective space. We identify explicit geometric generators: the Oakley--Usher torus equipped with three critical local systems, together with the Lagrangian sphere $\Sigma$.

\subsection{Outline and conventions}
Section~\ref{sec:tori_quadric} develops the~$\Psi$-invariant, constructs the Gelfand--Cetlin fibration, computes the two OU~$\Psi$-germs, and compares~$\cp^{\Ch}$ with lifted Vianna tori, proving Theorem~\ref{thm:quadric_exotic}. Section~\ref{sec:enumerative} computes the superpotential, proves non-displaceability and split-generation, and studies the change under compactification, proving Theorem~\ref{thm:quadric_non_displ}.

We use~$S^1=\R/\Z$ and integer coefficients for homology unless otherwise specified. Symplectic reduction is denoted by~$X\sslash_{\mu=c}G$. Induced moment maps on reduced or cut spaces retain the original notation. The symplectic normalization and size parameters agree with those in~\cite{OU1}; the Novikov field and local-system conventions are specified in Section~\ref{sec:enumerative}.

\subsection*{Acknowledgments} 
This paper and its companion~\cite{OU1} grew out of a Master's thesis written by AP under the supervision of JB at ETH Zürich. We thank Paul Biran for orchestrating this. We warmly thank Leonid Polterovich for drawing our attention to Oakley--Usher Lagrangians, and Johanna Bimmermann, Marc Fares, Johannes Hauber, David Keren Yaar, Joel Schmitz, Mike Usher, and Renato Vianna for interesting discussions, comments and questions. We are particularly grateful to Felix Schlenk for asking many interesting questions. JB is supported by SNSF Ambizione Grant PZ00P2-223460.

\subsection*{Statement on the use of LLMs} An earlier version of this project comprising this paper and its companion \cite{OU1} was written completely without the use of large language models. In subsequent versions, GPT was used to help split the project into two papers, and proofread the text. The authors take full responsibility for the final mathematical content.

\section{Chekanov OU-torus in the three-dimensional quadric}
\label{sec:tori_quadric}

In this section, we exclusively work with the case~$k=m=1$ in the three-dimensional quadric~$Q^3$. For ease of notation, we drop the~$1,1$-subscripts and sizes of the OU-Lagrangians everywhere. In other words, we write~$\cp = \cp_{1,1}(\frac{4\pi}{3})$ and~$\cp^{\Ch} = \cp_{1,1}^{\Ch}(\frac{4\pi}{3})$. Recall that these are the monotone OU-Lagrangian and Chekanov OU-Lagrangian, respectively.

\subsection{The \texorpdfstring{$\Psi$}{Psi}-invariant and its versal deformation} In order to distinguish tori up to Hamiltonian isotopy, we use the so-called $\Psi$-invariant, introduced by Shelukhin--Tonkonog--Vianna~\cite{SheTonVia24}. Let $(X,\omega)$ be a symplectic manifold.
For the general statements in this subsection, take $X$ closed and $L$ closed, orientable and spin. These hypotheses suffice for all applications here; extensions to open manifolds require the appropriate compactness hypotheses for the Fukaya algebra.

 The $\Psi$-invariant associates to every orientable spin Lagrangian~$L \subset X$ a number~$\Psi(L) \in (0,+\infty]$ which measures the smallest symplectic area of a~$J$-holomorphic disk which contributes nontrivially to some algebraic disk-count in the Fukaya~$A_{\infty}$-algebra of~$L$.

As we will see, $\Psi$ itself does not distinguish the tori we are interested in. For example, we will see that $\Psi(\cp) = \Psi(\cp^{\rm Ch}) = \frac{2\pi}{3}$. Therefore, we consider its versal deformation. We refer to \cite[Section 3]{BreHauSch23} for details on versal deformations. 

\begin{definition}
	Let~$L \subset X$ be an orientable spin Lagrangian. We call \emph{$\Psi$-germ (of $L$)} the germ, in the sense of \cite[Definition 3.13]{BreHauSch23}, of the invariant~$\Psi$ at~$L$. We denote it by
    	\begin{equation*}
        	\psi_L = [\Psi]_L \colon H^1(L;\R) \dashrightarrow \R \cup \{+\infty\},
    	\end{equation*}
    where the dashed arrow denotes a map defined on a neighbourhood of the origin of~$H^1(L;\R)$.
\end{definition}

\begin{proposition}
	\label{prop:min_functionals}
	The $\Psi$-germ is of the form 
	\begin{equation}
		\psi_L(a) = \Psi(L) + \min\{\langle a , \xi_i \rangle \},
	\end{equation}
	for some $\xi_1,\ldots,\xi_N \in H_1(L)$, provided~$\Psi(L)<+\infty$. In particular, the displayed piecewise-affine expression extends to $\psi_L \colon H^1(L;\R) \rightarrow \mathbb R$.
\end{proposition}

\proof
It was proved in \cite{SheTonVia24} that, in a small neighbourhood of $L$ in the space of Lagrangians, the invariant $\Psi$ can be written as the minimum of symplectic areas of relative classes $\beta_1,\ldots,\beta_N \in H_2(X,L)$ with boundary on $L$. 
More precisely, the classes $\beta_i$ are those of area
$\Psi(L)$ which contribute nontrivially to some
symmetrized $A_\infty$-operation on odd-degree
cohomology classes, in a classically minimal model
of the Fukaya algebra of $L$; see
\cite[Theorem 4.8 and the proof of Theorem 4.9]{SheTonVia24}.
 By Gromov's compactness theorem, this set is finite. Setting $\xi_i = \pp \beta_i$ and computing the versal deformation yields the claim. 
\proofend

Instead of working with $\psi_L$ directly, it is often more convenient to work with its (super-) level sets, since these are easier to manipulate and can be represented visually.

\begin{definition}
Let $L\subset X$ be an orientable spin Lagrangian
with $\Psi(L)<+\infty$.
Using the global piecewise-affine extension of
$\psi_L$ given by Proposition~\ref{prop:min_functionals},
we define
\begin{equation*}
    \Delta^{\!\Psi}(L)
    :=
    \{a\in H^1(L;\R)\mid \psi_L(a)\geqslant 0\}.
\end{equation*}
This is a convex polyhedron, which we call the
\emph{$\Psi$-polytope of $L$} whenever it is bounded. If $\Psi(L)=+\infty$, the germ is identically $+\infty$ by \cite[Theorem 4.9]{SheTonVia24}; one may separately set $\Delta^{\!\Psi}(L)=H^1(L;\R)$.

\end{definition}

\begin{remark}
	For the Oakley--Usher torus considered in the next subsection, the relevant structure is a \emph{Gelfand--Cetlin fibration}, meaning that the $\Psi$-polytope agrees with 
	\begin{enumerate}
		\item the polytope obtained as the base of the fibration, after a translation that shifts the image of the monotone fibre $L$ to the origin, if necessary.
		\item the appropriately scaled polar dual of the Newton polytope of the superpotential of that fibre.
	\end{enumerate}
	The fact that these three polytopes agree is \cite[Theorem B]{SheTonVia24}. We only use Gelfand--Cetlin fibrations once to compute $\psi_\cp$ for the Oakley--Usher Lagrangian $\cp \subset Q^3$. All other invariants will be computed by versal deformation techniques. 
\end{remark}

We use the following result to distinguish tori.

\begin{proposition}
	\label{prop:distinguish_tori}
	Let $L,L' \subset X$ be orientable spin Lagrangian manifolds. 
	If a symplectomorphism~$\phi:X\to X$ maps~$L$ to~$L'$, then there exists an isomorphism $\chi \colon H^1(L;\Z) \rightarrow H^1(L';\Z)$ such that $\chi_\R(\Delta^{\!\Psi}(L)) = \Delta^{\!\Psi}(L')$, where $\chi_\R \colon H^1(L;\R) \rightarrow H^1(L';\R)$ denotes the extension to $\R$-coefficients.
\end{proposition}

\proof
Set~$\chi=(\phi|_L^*)^{-1}$. Pushing forward the almost complex structure and transporting the spin structure by~$\phi$ identifies the relevant disk moduli spaces and their symplectic areas, so the definition of~$\Psi$ in~\cite{SheTonVia24} gives~$\Psi(\phi(K))=\Psi(K)$. The flux of a deformation is carried to its image under~$\chi_\R$; consequently~$\psi_{L'}=\psi_L\circ\chi_\R^{-1}$. Taking nonnegative superlevel sets gives the claim.
\proofend

In practice, we usually identify $H^1(L;\R)$ and $H^1(L';\R)$ with $\R^n$ and consider the polytopes up to $\GL(n;\Z)$. In that case, we obtain that the $\Psi$-polytopes of $L$ and $L'$ are related by an element in $\GL(n;\Z)$ whenever an ambient symplectomorphism takes one to the other.

\subsection{The Oakley--Usher Lagrangian}
To compute the $\Psi$-germ of $\cp$, we rely on \cite[Proposition C]{SheTonVia24}. It determines $\Psi$ of fibres of so-called \emph{Gelfand--Cetlin fibrations} which have a monotone fibre. 

\begin{figure}[h!]
\centering
\begin{tikzpicture}[scale=1.05, >=Latex, line cap=round, line join=round]

\coordinate (O) at (0,0);
\coordinate (e1) at (2.5,0.6);
\coordinate (e2) at ( 2,-0.7);
\coordinate (e3) at ( 0.0,3);

\draw[->, thick, gray!70] (O) -- ($(O)+1.2*(e1)$) node[below right] {$x_2$};
\draw[->, thick, gray!70] (O) -- ($(O)+1.2*(e2)$) node[below] {$x_1$};
\draw[->, thick, gray!70] (O) -- ($(O)+1.2*(e3)$) node[above] {$x_3$};
\coordinate (T) at ($(O)+(e3)$);
\coordinate (A) at ($(O)+(e1)+(e3)$);
\coordinate (B) at ($(O)+(e2)+(e3)$);
\coordinate (C) at ($(O)+(e1)+(e2)+(e3)$);
\coordinate (L) at ($(O)+0.9*(e1)+0.9*(e2)+0.7*(e3)$);

\coordinate (P) at ($(O)+0.303*(e1)+0.333*(e2)+0.667*(e3)$);

\fill[gray!5 ]  (T)--(A)--(C)--(B)--cycle;
\fill[gray!20 ]  (O)--(T)--(B)--cycle;
\fill[gray!40] (O)--(B)--(C)--cycle;

\draw[thick] (T)--(A)--(C)--(B)--cycle;
\draw[thick] (O)--(T)--(B)--cycle;
\draw[thick] (O)--(B)--(C)--cycle;

\draw[dashed] (O)--(A);

\fill[blue!80!black] (P) circle (2.2pt);
\node[blue!80!black, left=3pt] at (P)
{$x_{\mathcal{P}}$};

\node[font=\small] at (L) {$\Delta_{Q^3}$};

\end{tikzpicture}

\caption{ \label{fig:GCpol} The Gelfand--Cetlin polytope~\(\Delta_{Q^3}\)}
\end{figure}
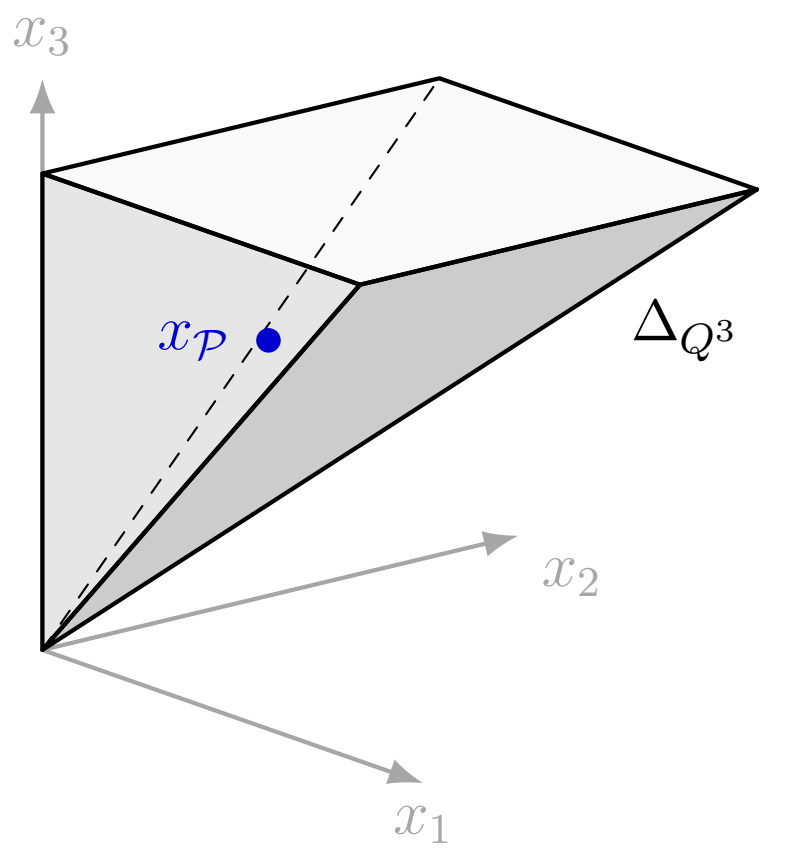

A continuous Lagrangian torus fibration is called \emph{Gelfand--Cetlin} if it restricts to a toric moment map away from the codimension two faces.
 As it turns out, the moment map~$\mu_{1,1}$ defined by~\cite[Proposition 3.2]{OU1} induces a Gelfand--Cetlin fibration on~$Q^3$ in which~$\cp$ appears as a monotone fibre, see \cref{fig:GCpol} for its polytope.

\begin{proposition}
    \label{prop:GCfibration}
    There is a Gelfand--Cetlin fibration~$\mu \colon Q^3 \rightarrow \Delta_{Q^3}$ over the polytope 
    \begin{equation}
        \label{eq:deltaQ3}
        \Delta_{Q^3} =
        \{x_1,x_2 \geqslant 0, \; x_3 - x_1 \geqslant 0, \; x_3 - x_2 \geqslant 0, \; 2\pi-x_3 \geqslant 0\} \subset \R^3,
    \end{equation}
    such that
    \begin{enumerate}
        \item the fibre~$\Sigma = \mu^{-1}(0)$ is a Lagrangian sphere;
        \item the fibre~$\cp = \mu^{-1}(x_{\cp})$ of~$x_{\cp} = (\frac{2\pi}{3},\frac{2\pi}{3},\frac{4\pi}{3})$ is the Oakley--Usher Lagrangian, and 
        \item the restriction~$\mu \colon Q^3 \setminus \Sigma \rightarrow \Delta_{Q^3} \setminus \{0\}$ is a smooth toric moment map.
    \end{enumerate}
\end{proposition}

\proof
Recall that the continuous moment map 
\begin{equation}
    \label{eq:mu11}
    \mu_{1,1} \colon T^*S^3 \rightarrow \R^3, \quad (q,p) \mapsto \bigl(2\pi \vert p \vert, 2\pi(p_1q_2 - p_2q_1), 2\pi(p_3q_4 - p_4q_3)\bigr)
\end{equation}
generates a Hamiltonian~$T^3$-action away from the zero-section~$0_{S^3}$. Here the last two scalar components use period-one generators of the two~$\SO(2)$ factors; this specifies their normalization and sign in the identification with~$\R^2$. This moment map is smooth away from the zero-section, and its image is the cone defined by the following four inequalities
    \begin{equation*}
        \widehat{\Delta}_{T^*S^3}
        =
        \{x_1 \pm x_2 \pm x_3 \geqslant 0\} \subset \R^3.
    \end{equation*}
Recall however from \cite[Remark 3.1]{OU1} that the element~$\theta_{1/2} = (1/2,1/2,1/2) \in T^3$ acts trivially, meaning that the action is not effective and thus not toric. This can be remedied as follows: Viewing the standard torus~$T^3$ as the quotient of~$\R^3$ by the lattice~$\Lambda_0 = \Z^3$, we can add~$\theta_{1/2}$ as a generator to the lattice~$\Lambda_0$ to obtain~$\Lambda$. The new torus~$\R^3 / \Lambda$ acts effectively and has the same image under~$\mu_{1,1}$, though with respect to a different lattice in~$\R^3 = \mathfrak{t}^*$ (the target space of~$\mu_{1,1}$). The dual lattice of~$\Lambda$ is given by
\begin{equation*}
    \Lambda^* := \{(x_1,x_2,x_3) \in \mathbb{Z}^3 \sth x_1 + x_2 + x_3 \in 2\Z \}.
\end{equation*}
Since we are used to drawing the image of the moment map in~$\R^3 = \mathfrak{t}^*$ equipped with the \emph{standard} integral affine structure~$\Lambda_0^* = \Z^3$, we apply some~$M \in \GL(3;\Q)$ which induces an isomorphism from~$\Lambda^*$ to~$\Z^3$. It is not hard to show that the linear transformation
\begin{equation}
    \label{eq:basechange_lattice}
    M = \frac{1}{2}
    \begin{pmatrix}
        1 & -1 & -1 \\
        1 & 1 & -1 \\
        2 & 0 & 0
    \end{pmatrix}
\end{equation}
has this property. We define a new moment map~$\mu = M \mu_{1,1}$. Since the kernel of the original torus action is exactly $\{0,\theta_{1/2}\}$, this moment map induces an effective Hamiltonian~$T^3$-action and turns~$T^*S^3 \setminus 0_{S^3}$ into a toric manifold. Its continuous extension to~$0_{S^3}$ satisfies~$\mu^{-1}(0) = 0_{S^3}$ and its image is the cone
    \begin{equation*}
        \Delta_{T^*S^3}
        =
        \{x_1,x_2 \geqslant 0, \; x_3 - x_1 \geqslant 0, \; x_3 - x_2 \geqslant 0\} 
        \subset \R^3.
    \end{equation*}
Recall from~\cite[Proposition~3.10]{OU1} that performing a symplectic cut at~$2\pi \vert p \vert = 2\pi$ yields the quadric~$Q^3$. The Hamiltonian~$2\pi \vert p \vert$ is equal to the first component of the original moment map~$\mu_{1,1}$ and thus corresponds to the third component of~$\mu$, as a computation using~$M$ shows. Therefore, performing a symplectic cut with respect to~$\mu_3 = 2\pi$ yields~$Q^3$. By abuse of notation, we denote the resulting moment map by~$\mu \colon Q^3 \rightarrow \R^3$, too. By construction, its image is~$\Delta_{Q^3}$ as defined in \eqref{eq:deltaQ3}. Let~$\Sigma = \mu^{-1}(0) \subset Q^3$. This is the Lagrangian sphere obtained from~$0_{S^3} \subset T^*S^3$ after the symplectic cut. 

Recall from~\cite[Proposition~3.4, equation~(17)]{OU1} that the Oakley--Usher Lagrangians~$\cp_{1,1}(r) \subset T^*S^3$ are given by~$\mu_{1,1}^{-1}(r,0,0)$ in the original moment map and thus we obtain 
\begin{equation*}
    \cp_{1,1}(r) = \mu^{-1}(\textstyle\frac{r}{2},\textstyle\frac{r}{2},r),
\end{equation*}
with respect to the new moment maps for all~$r > 0$ in~$T^*S^3$ and for~$r \in (0,2\pi)$ in~$Q^3$. Recall from~\cite[Proposition~3.16]{OU1} that the Oakley--Usher Lagrangian in~$Q^3$ is monotone for~$r = \frac{4\pi}{3}$. It follows that~$\cp = \mu^{-1}(x_{\cp}) \subset Q^3$ for~$x_{\cp}=(\frac{2\pi}{3},\frac{2\pi}{3},\frac{4\pi}{3})$. By construction, the moment map~$\mu \colon Q^3 \rightarrow \Delta_{Q^3}$ is continuous and a smooth toric moment map away from~$\Sigma$.
\proofend

Let 
    \begin{equation*}
        T(x) = \mu^{-1}(x), \quad
        x \in \Int \Delta_{Q^3},
    \end{equation*}
denote the corresponding Gelfand--Cetlin fibres and let~$d_{\rm IA}(x,\pp \Delta_{Q^3})$ be the integral affine distance of the point~$x \in \R^3$ to the boundary of~$\Delta_{Q^3}$. Explicitly, it is given by
    \begin{equation}
    	\label{eq:IAdist_Q^3}
        d_{\rm IA}(x,\pp \Delta_{Q^3})
        =
        \min \{x_1,x_2,x_3 - x_1, x_3 - x_2, 2\pi-x_3\}.
    \end{equation}

Together with \cref{prop:GCfibration}, we use \cite[Proposition C]{SheTonVia24} to find the following. 

\begin{proposition}
    \label{prop:Psi_fibres}
    For every~$x \in \Int \Delta_{Q^3}$, the~$\Psi$-invariant of a Gelfand--Cetlin fibre~$T(x)$ is equal to the integral affine distance of~$x$ to the boundary~$\pp \Delta_{Q^3}$,
    \begin{equation*}
        \Psi(T(x)) = d_{\rm IA}(x,\pp \Delta_{Q^3}).
    \end{equation*}
\end{proposition}

\begin{remark}
    The ease with which the results in \cite{SheTonVia24} allow us to compute~$\Psi$ on all of the fibres of~$\mu$ is the reason why we use the invariant~$\Psi$ instead of displacement energy. In fact, we suspect that whenever~$\Psi$ is affine linear (as opposed to just piecewise affine linear) on a neighbourhood of a point~$x \in \Int \Delta_{Q^3}$, then the displacement energy agrees with~$\Psi$. The upper bound on displacement energy readily follows from probes (see for example \cite{McD11}) in the Gelfand--Cetlin polytope, but the lower bound requires more thought.
\end{remark}

In particular, \cref{prop:Psi_fibres} allows us to compute the~$\Psi$-germ of~$\cp$. Let~$H_1(\cp) \cong \Z^3$ be the isomorphism induced by taking the generators of~$H_1(\cp)$ defined by the (local) toric action generated by~$\mu$. By \cite[Proposition 3.15]{BreHauSch23}, we obtain a versal deformation of~$\cp$ by considering~$a \mapsto \mu^{-1}(x_{\cp} + a)$, where~$a \in \R^3 \cong H^1(\cp;\R)$ under the identification dual to the above isomorphism. Recall that $\cp = \mu^{-1}((\frac{2\pi}{3},\frac{2\pi}{3},\frac{4\pi}{3}))$ in the Gelfand--Cetlin fibration. From Proposition \ref{prop:Psi_fibres} and \eqref{eq:IAdist_Q^3} we obtain the following.

\begin{corollary}
\label{cor:germP}
There is an identification~$H^1(\cp;\R) \cong \R^3$ such that for all $a \in \R^3$,
    \begin{align*}
    	\psi_{\cp}(a) &= \textstyle\frac{2\pi}{3} + \min\{a_1,a_2,a_3 - a_1, a_3 - a_2,- a_3\},\\
    	\Delta^{\!\Psi}(\cp) &= \{a \in \R^3 \sth a_1,a_2,a_3-a_1,a_3-a_2,-a_3 \geqslant -\textstyle\frac{2\pi}{3} \}.
    \end{align*}
\end{corollary}

As expected, the $\Psi$-polytope of $\cp$ agrees with the base polytope $\Delta_{Q^3}$ of the Gelfand--Cetlin fibration after a translation by $-\mu(\mathcal{P})=(-\frac{2\pi}{3},-\frac{2\pi}{3},-\frac{4\pi}{3})$.

\subsection{The Chekanov Oakley--Usher Lagrangian}
\label{ssec:pchek}
On the other hand, we obtain the following~$\Psi$-germ for the Chekanov Oakley--Usher Lagrangian~$\cp^{\rm Ch}$. 

\begin{proposition}
    \label{prop:germPchek}
    There is an identification~$H^1(\cp^{\rm Ch};\R) \cong \R^3$ such that for all $b \in \R^3$,
    \begin{align*}
        \psi_{\cp^{\Ch}}(b) &= \textstyle\frac{2\pi}{3} + \min\{\pm b_1 \pm b_2 - 2b_3, b_3\},\\
        \Delta^{\!\Psi}(\cp^{\Ch}) &= \{b \in \R^3 \sth \pm b_1 \pm b_2 - 2b_3, b_3 \geqslant -\textstyle\frac{2\pi}{3}\}.
    \end{align*}
\end{proposition}

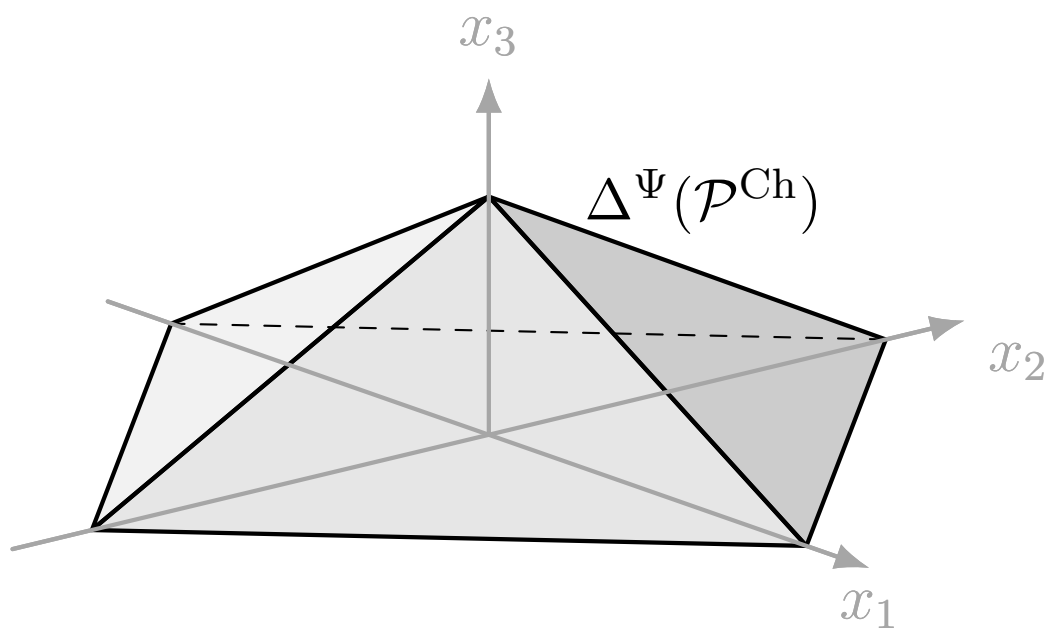
\begin{figure}[h!]
\centering
\begin{tikzpicture}[scale=1.05, >=Latex, line cap=round, line join=round]

\coordinate (O) at (0,0);
\coordinate (e1) at (2.5,0.6);
\coordinate (e2) at ( 2,-0.7);
\coordinate (e3) at ( 0.0,3);

\draw[->, thick, gray!70] ($(O)-1.2*(e1)$) -- ($(O)+1.2*(e1)$) node[below right] {$x_2$};
\draw[->, thick, gray!70] ($(O)-1.2*(e2)$) -- ($(O)+1.2*(e2)$) node[below] {$x_1$};
\draw[->, thick, gray!70] (O) -- ($(O)+0.75*(e3)$) node[above] {$x_3$};
\coordinate (T) at ($(O)+0.5*(e3)$);
\coordinate (A) at ($(O)-(e1)$);
\coordinate (B) at ($(O)+(e2)$);
\coordinate (C) at ($(O)+(e1)$);
\coordinate (D) at ($(O)-(e2)$);
\coordinate (L) at ($(O)+0.3*(e1)+0.3*(e2)+0.5*(e3)$);

\fill[gray!10 ]  (A)--(T)--(D)--cycle;
\fill[gray!20 ]  (A)--(B)--(T)--cycle;
\fill[gray!40] (B)--(C)--(T)--cycle;

\draw[thick] (A)--(T)--(D)--cycle;
\draw[thick] (A)--(B)--(T)--cycle;
\draw[thick] (B)--(C)--(T)--cycle;

\draw[dashed] (C)--(D);

\draw[->, thick, gray!70] ($(O)-1.2*(e1)$) -- ($(O)+1.2*(e1)$) node[below right] {$x_2$};
\draw[->, thick, gray!70] ($(O)-1.2*(e2)$) -- ($(O)+1.2*(e2)$) node[below] {$x_1$};
\draw[->, thick, gray!70] (O) -- ($(O)+0.75*(e3)$) node[above] {$x_3$};

\node[font=\small] at (L) {$\Delta^{\Psi}(\cp^{\rm Ch})$};

\end{tikzpicture}

\caption{ \label{fig:Psipol} The $\Psi$-polytope of the Chekanov OU Lagrangian $\cp^{\rm Ch}$. We have applied a translation in the $x_3$-direction for readability.}
\end{figure}

See \cref{fig:Psipol} for an illustration of $\Delta^\Psi(\cp^{\rm Ch})$. The $\Psi$-polytopes $\Delta^{\!\Psi}(\cp)$ and $\Delta^{\!\Psi}(\cp^{\Ch})$ cannot be mapped to one another by an element in $\GL(3;\Z)$. One way of seeing this is as follows: The integral affine distance of the four-valent vertex of $\Delta^{\!\Psi}(\cp)$ to its opposite facet is $2\pi$. For $\Delta^{\!\Psi}(\cp^{\Ch})$ this distance is $\pi$. This is well-defined and an invariant under $\GL(n;\Z)$. By Proposition \ref{prop:distinguish_tori}, we obtain the following.

\begin{theorem}
	No ambient symplectomorphism of $Q^3$ maps $\cp$ to $\cp^{\Ch}$.
\end{theorem}

All that is left to prove is \cref{prop:germPchek}. Its proof follows the same pattern as the computation of displacement energy germs in \cite{Bre22, Bre23b, BreHauSch23, CheSch10}: on a small neighbourhood of~$\cp^{\Ch}$ in its versal deformation, there is an open dense set on which a deformation of~$\cp^{\Ch}$ is Hamiltonian isotopic to a deformation of~$\cp$. Note that the latter are just the fibres $T(x)$ of the Gelfand--Cetlin fibration and thus we can use \cref{prop:Psi_fibres} to prove \cref{prop:germPchek}. This computes $\psi_{\cp^{\Ch}}$ and the $\Psi$-polytope on an open dense subset of a neighbourhood of the origin. By Proposition~\ref{prop:min_functionals}, this determines $\psi_{\cp^{\Ch}}$ and the polytope everywhere. \smallskip

\proofof{\cref{prop:germPchek}} The proof is very similar to Sections 4 and 5 in \cite{Bre23b}, which we recommend for more details. 
\\
\\
\noindent
{\bf Step 0: Outline.} Recall that, by construction,
the Lagrangian torus~$\cp^{\Ch}$ is the lift of a
curve~$\Gamma(\frac{4\pi}{3})$ bounding area~$\frac{2\pi}{3}$
in the smooth part of the symplectic reduction
\[
Q^3\sslash_{\overline{\mu}_{1,1}=0}G_{1,1},
\]
see the discussion surrounding~\cite[Definitions~3.15 and~3.20]{OU1}.
Strictly speaking, this reduction differs from the
space~$\widehat Z$ considered in~\cite{OU1}:
the actual quotient contains the interval
$\Sigma/G_{1,1}\cong[0,1]$, whereas
$\widehat Z$ is obtained by collapsing this interval
to a point. Since the curve~$\Gamma(\frac{4\pi}{3})$
and its enclosed disk lie away from this interval,
the two descriptions agree on the region relevant
to our construction.

Recall that a versal deformation of the torus~$\cp^{\Ch} \subset Q^3$ is a continuous map 
\begin{equation}
    \label{eq:vdcpchek}
    v_{\cp^{\Ch}} \colon H^1(\cp^{\Ch};\R) \cong \R^3 \dashrightarrow \cl_{\cp^{\Ch}}
\end{equation}
defined on a neighbourhood of the origin~$0 \in H^1(\cp^{\Ch};\R)$ which is a local section to the Lagrangian flux map. By~$\cl_{\cp^{\Ch}}$ we have denoted the space of Lagrangian deformations of~$\cp^{\Ch}$ equipped with the~$C^1$ topology. For a detailed exposition, we recommend \cite[\S 3.3]{BreHauSch23}, see in particular Definition 3.9 therein. To define such a versal deformation of~$\cp^{\Ch}$ explicitly, we will carry out symplectic reduction on a \emph{family of probes} pointing in the direction~$(1,1,2)$ and identify the members of the versal deformation as certain lifts of curves in the family of reduced spaces. The deformation space of~$\cp^{\Ch}$ is three-dimensional, see \eqref{eq:vdcpchek}: one of the parameters will come from varying the area of the curve~$\Gamma(\cdot)$ we lift, and the two other parameters will come from varying the reduced space in the family that we lift the curve from. \smallskip
\noindent

{\bf Step 1: Symplectic reduction.} Let us now set up the family of symplectic reductions. To that end, we define the following moment map
\begin{equation}
    \label{eq:nu}
    \nu \colon Q^3 \rightarrow \R^2, \quad
    \nu = (\mu_2 - \mu_1, \mu_3 - \mu_1 - \mu_2 ).
\end{equation}
Let us point out that the level sets (fibres) of~$\nu$ point in the direction~$(1,1,2)$ in the image of~$\mu$, which justifies this choice. Since its components are integral linear combinations of the components of~$\mu$, the moment map~$\nu$ defines a~$T^2$-subaction of the~$T^3$-action generated by~$\mu$. In fact, the inclusion of tori~$T^2 \hookrightarrow T^3$ is given by~$(\theta_1,\theta_2) \mapsto (-\theta_1 -\theta_2 , \theta_1 -\theta_2 , \theta_2)$. Let us compute the image of~$\nu$. To that end, recall from \cref{prop:GCfibration} that the image~$\Delta_{Q^3}$ of~$\mu$ is defined by a set of five inequalities~$\ell_i(x) \geqslant 0$ for 
\begin{equation}
    \label{eq:ell}
    \ell_0(x) = 2\pi - x_3, \quad
    \ell_1(x) = x_1,\quad
    \ell_2(x) = x_2, \quad
    \ell_3(x) = x_3 - x_1\quad
    \ell_4(x) = x_3 - x_2. 
\end{equation}
A computation shows that under these constraints the image of~$\nu = (x_2-x_1,x_3-x_1-x_2)$ is given by the lozenge
\begin{equation}
    \Diamond = \{(b_1,b_2) \in \R^2 \sth \pm b_1 \pm b_2 \leqslant 2\pi \}.
\end{equation}
Indeed, since~$\Diamond = \im (\nu)$ is the image of a linear projection of the convex polytope~$\Delta_{Q^3}$, it is itself a convex polytope spanned by the image of the vertices of~$\Delta_{Q^3}$. Now let~$(b_1,b_2) \in \Int (\Diamond)$ and let 
\begin{equation}
    \label{eq:nureduction}
    \pi_{(b_1,b_2)} \colon \nu^{-1}(b_1,b_2) \rightarrow S_{(b_1,b_2)}
\end{equation}
be the corresponding, possibly singular,
symplectic reduction by the~$T^2$-action
generated by~$\nu$. We think of this as the reduction associated to the \emph{singular probe}
\begin{equation}
\label{eq:singular probe}
    P_{(b_1,b_2)} = \{(x_2 - x_1, x_3 - x_1 - x_2) = (b_1,b_2)\} \cap \Delta_{Q^3}.
\end{equation}
Recall that the intersection of a probe with the boundary of a toric moment polytope yields a smooth point in the reduced space if and only if this intersection happens in the interior of a facet and satisfies the \emph{integral transversality} condition. This condition asks that a lattice basis in the facet in question can be completed to a lattice basis of the full ambient lattice by the direction vector of the probe, see also \cite[\S 2.1]{McD11} or \cite[Definition 2.3]{Bre23} for a generalization. The singular probes~$P_{(b_1,b_2)}$ turn out to be honest probes for all~$(b_1,b_2)$ with~$b_1,b_2 \neq 0$. This follows from the fact that their direction vector~$(1,1,2)$ intersects all four facets~$F_i = \{\ell_i = 0\}$ for~$i \in \{1,2,3,4\}$ integrally transversely. Here~$\ell_i$ is defined as in~\eqref{eq:ell}. Its intersection with~$F_0 = \{ \ell_0 = 0\}$ is not integrally transverse, leading to a~$\Z_2$-orbifold point in all reduced spaces (as in the case of~$\widehat{Z}$ in general). Let us discuss what the reduced spaces look like in more detail. 
\begin{enumerate}

\item In the case $b_1=b_2=0$, we recover
the reduction used in the construction of
$\cp^{\Ch}$, with one distinction.
The actual quotient $S_{(0,0)}$ contains
the interval $\Sigma/T^2\cong[0,1]$, parametrized by $|q_1|^2+|q_2|^2$.
The space $\widehat Z$ considered in
\cite{OU1} is obtained by collapsing
this interval to a point. It is a
two-sphere with a $\Z_2$-orbifold point
at the upper end and a distinguished
point arising from this collapse.
Since the curve defining $\cp^{\Ch}$
and its enclosed disk lie in the smooth
part away from the interval, the two
descriptions agree on the region
relevant to our construction.

\item In case only one of the~$b_i$ is zero, the reduced space is a two-sphere containing a~$\Z_2$-orbifold point and a distinguished point coming from the fact that the probe~$P_{(b_1,b_2)}$ intersects one of the (toric!) edges of~$\Delta_{Q^3}$, whose fibre is a circle. Its smooth part is symplectomorphic to a standard cylinder. 
\item In case~$b_1b_2 \neq 0$, the reduced space is a two-sphere containing only the~$\Z_2$-orbifold point coming from the upper facet. The complement of that point is symplectomorphic to a standard symplectic two-disk. Again, the smoothness in the interior of the disk is due to the fact that the probe~$P_{(b_1,b_2)}$ intersects one of the facets~$F_1,\ldots,F_4$ integrally transversely.
\end{enumerate}
\medskip

\noindent
{\bf Step 2: Versal deformation.} Let us now define the versal deformation of~$\cp^{\Ch}$. We recall (for example from \cite[\S 3.3]{BreHauSch23}) that any Weinstein chart~$\varphi \colon T^*L \supset U \rightarrow X$ of a Lagrangian torus~$L \subset X$ induces a versal deformation~$v_L^{\varphi}(b) = \varphi(\{b\} \times T^n)$. Here we have identified~$T^*L = \R^n \times T^n$ and~$H^1(L;\R) = \R^n$. To make our versal deformation compatible with the symplectic reduction set up above, we choose a~$T^2$-equivariant Weinstein chart of~$\cp^{\Ch}$. On the model space~$T^*T^3 = \{(p_i,q_i)\}$, we pick the~$T^2$-action defined by standard rotation in the~$q_1,q_2$ components generated by the moment map given by the projection to~$(p_1,p_2)$. On~$Q^3$, we have the~$T^2$-action defined by~$\nu$ as in \eqref{eq:nu}.
Although $\mu$ is not smooth at $\Sigma$, the two components of $\nu$ are the last two components of $\mu_{1,1}$. Hence $\nu$ is smooth everywhere and generates the original $T^2$-action.

 Since~$\cp^{\Ch}$ is invariant under the latter, we can choose a Weinstein chart~$\varphi$ which is~$T^2$-equivariant with respect to these two actions. Define a versal deformation 
\begin{equation}
    \label{eq:VDcpchekdef}
    v = v_{\cp^{\Ch}}^{\varphi} \colon H^1(T^3;\R) \cong \R^3 \dashrightarrow \cl_{\cp^{\Ch}} , \quad
    v(b) = \varphi(\{b\} \times T^3).
\end{equation}

\begin{lemma}
    \label{lem:VDlift}
    The Lagrangians appearing in the versal deformation \eqref{eq:VDcpchekdef} are lifts of curves from the family of symplectic reductions \eqref{eq:nureduction}. More specifically, there exists a family of curves~$\mathcal{C}_{(b_1,b_2)}(b_3) \subset S_{(b_1,b_2)}$ contained in the smooth part of the reduced spaces~$S_{(b_1,b_2)}$ and bounding a disk of symplectic area~$\frac{2\pi}{3} + b_3$ such that 
    \begin{equation*}
        \label{eq:VDlift}
        v(b) = \pi^{-1}_{(b_1,b_2)}(\mathcal{C}_{(b_1,b_2)}(b_3))
    \end{equation*}
    for all~$b = (b_1,b_2,b_3) \in H^1(T^3;\R) \cong \R^3$ contained in a small enough neighbourhood $V$ of the origin such that~$v$ is defined.
\end{lemma}

\proof
By equivariance of $\phi$, the tori $v(b)$ are
$T^2$-invariant with respect to the action generated
by $\nu$. Therefore, these tori are lifts of curves
in the corresponding reduced spaces.

To compute the areas of the disks bounded by these
curves, we choose the first two generators of
$H_1(\mathcal P^{\mathrm{Ch}};\Z)$ to be the
$T^2$-orbit classes. We choose the third generator
to be the boundary of a section over the disk
bounded by the original reduced curve. Such a
section exists since the principal $T^2$-bundle
over this disk is trivial. With this choice of marking, the equivariant
Weinstein chart can be normalized so that
$\nu\circ\phi=(p_1,p_2)$. For sufficiently small
$b$, the reduced curves and their bounding disks
remain in the smooth part of the reduced spaces,
and the corresponding section disks can be chosen
to vary smoothly.

Their symplectic areas agree with the reduced
disk areas, and their variation is given by the
Lagrangian flux evaluated on the third generator.
By the definition of the versal deformation, this
variation is $b_3$. Since the original reduced disk
has area $2\pi/3$, we obtain
\[
\operatorname{Area}(D_b)=\frac{2\pi}{3}+b_3.
\]
The argument is the same as in
\cite[Proposition 4.4]{Bre23b}, with the above
choice of marking and continuously varying
section disks.
\proofend

\noindent
{\bf Step 3: Deforming the tori.} The crucial observation which makes this proof work is: whenever~$b_1,b_2 \neq 0$, the reduced space~$S_{(b_1,b_2)}$ is a sphere containing a single~$\Z_2$-orbifold point (and no distinguished points), see point (3) in the enumeration above. Therefore the curve~$\mathcal{C}_{(b_1,b_2)}(b_3) \subset S_{(b_1,b_2)}$ can be deformed by a Hamiltonian isotopy into a \emph{standard} curve in~$S_{(b_1,b_2)}$ bounding the same area. By \emph{standard} we mean an orbit of the residual~$S^1$-action on the reduced space, or, equivalently, the curve obtained as the image of the corresponding fibre~$T(x) = \mu^{-1}(x)$ under the reduction. This can be used to prove that~$v(b)$ is Hamiltonian isotopic to some~$T(x)$, and thus, by \cref{prop:Psi_fibres}, allows us to compute the~$\Psi$-invariant~$\Psi(v(b))$ whenever~$b_1,b_2 \neq 0$. More precisely, we will show the following. 

\begin{lemma}
\label{lem:bigPhi}
    We have
    \begin{equation*}
        v(b) \cong T(\Phi(b)) \quad
        \text{for all } b = (b_1,b_2,b_3) \text{ with } b_1,b_2 \neq 0 \text{ in } V
    \end{equation*}
    where~$\Phi \colon \R^3 \rightarrow \R^3$ is the continuous and bijective piecewise integral affine map defined by
    \begin{equation}
        \label{eq:bigPhi}
        \Phi(b_1,b_2,b_3) - \left( \textstyle\frac{2\pi}{3} , \textstyle\frac{2\pi}{3}, \textstyle\frac{4\pi}{3} \right)
        =
        \begin{cases}
            (b_3, b_1+b_3, b_1 +b_2 + 2b_3) & b_1 \geqslant 0, b_2 \geqslant 0 \\
            (b_3-b_1, b_3, b_2 - b_1 + 2b_3) & b_1 \leqslant 0, b_2 \geqslant 0 \\
            (b_3 -b_1 - b_2, b_3 - b_2, 2b_3 - b_1 - b_2) & b_1 \leqslant 0, b_2 \leqslant 0 \\
            (b_3 - b_2, b_1 - b_2 + b_3, b_1 - b_2 + 2b_3) & b_1 \geqslant 0, b_2 \leqslant 0.
        \end{cases}
    \end{equation}
\end{lemma}

\proof 
Note that the condition~$b_1,b_2 \neq 0$ separates~$\Int (\Diamond)$ into four disjoint open regions. Let us focus on the case~$b_1,b_2 > 0$. Since~$\nu = (x_2 - x_1, x_3 - x_1 -x_2)$, this implies~$x_2 > x_1$ and~$x_3 > x_1 + x_2$ in~$\Delta_{Q^3}$, meaning that this yields only probes~$P_{b_1,b_2}$ intersecting the facet~$F_1 = \{\ell_1 = 0\}$ on the lower and the facet~$F_0 = \{\ell_0 = 0\}$ on the upper end. As we have discussed, the lower intersection with~$\Delta_{Q^3}$ is integrally transverse whereas the upper one contributes a~$\Z_2$-orbifold point to the reduced space~$S_{(b_1,b_2)}$. By \cref{lem:VDlift}, we find that~$v(b)$ (for~$b_1,b_2 > 0$ and arbitrary~$b_3$ such that~$(b_1,b_2,b_3) \in V$) projects to a curve~$\mathcal{C}_{(b_1,b_2)}(b_3)$ bounding area~$\frac{2\pi}{3} + b_3$ in the reduced space~$S_{(b_1,b_2)}$. By a Hamiltonian isotopy we can map the curve~$\mathcal{C}_{(b_1,b_2)}(b_3)$ to the standard circle~$S^1(\frac{2\pi}{3} + b_3)$ bounding the same area in the disk obtained as the smooth part of~$S_{(b_1,b_2)}$. After reduction by~$T^2$, the~$T^3$-action on~$Q^3$ (we are away from~$\Sigma$) induces a residual~$S^1$-action on~$S_{(b_1,b_2)}$, which allows us to identify the smooth part of the reduced space~$S^1$-equivariantly with the standard symplectic disk carrying the standard~$S^1$-action by rotation. The standard circle~$S^1(\frac{2\pi}{3} + b_3)$ is an orbit of that residual~$S^1$-action. Its preimage~$\pi^{-1}_{(b_1,b_2)}(S^1(\frac{2\pi}{3} + b_3))$ under the reduction is a fibre torus~$T(x) = \mu^{-1}(x)$. By lifting this Hamiltonian isotopy, we find that~$v(b)$ is Hamiltonian isotopic to $T(x)$. 

Let us determine~$x = (x_1,x_2,x_3) =: \Phi(b)$ in terms of~$b$. Since by construction~$T(x)$ is contained in the level set~$\nu^{-1}(b_1,b_2)$, we find that~$x_2 - x_1 = b_1$ and~$x_3 - x_1 - x_2 = b_2$. This constrains~$x$ to lie on the probe~$P_{(b_1,b_2)}$ defined in \cref{eq:singular probe}. Since $T(x)$ projects to the circle~$S^1(\frac{2\pi}{3} + b_3)$ in the reduced space, the integral affine distance of~$x$ to the facet~$F_1 = \{\ell_1 = 0\}$ is~$\frac{2\pi}{3} + b_3$. Indeed, on the smooth part of the reduced space $S_{(b_1,b_2)}$, the residual $S^1$-action is Hamiltonian and identifies this space $S^1$-equivariantly with the standard symplectic disk of area equal to the affine length of the probe $P_{(b_1,b_2)}$. Under this identification, the moment map of the residual $S^1$-action is exactly the integral affine coordinate along the probe, normalized to vanish at the facet where the probe enters the polytope. Combining these facts, we find
\begin{equation}
    \Phi(b_1,b_2,b_3) 
    = 
    \left( \textstyle\frac{2\pi}{3} , \textstyle\frac{2\pi}{3}, \textstyle\frac{4\pi}{3} \right)
    + (b_3, b_1 + b_3, b_1 + b_2 + 2b_3). 
\end{equation}
The same kind of argument shows the formula in the three remaining cases. The map defined in \eqref{eq:bigPhi} is the unique continuous extension to the region~$b_1b_2 = 0$ of these partially defined maps.
\proofend

With this in hand, we can finish the proof of \cref{prop:germPchek}. Combining \cref{lem:bigPhi} and \cref{prop:Psi_fibres}, we find
\begin{equation}
    \label{eq:PSIcomp}
    \Psi(v(b))
    = \Psi(T(\Phi(b)))
    = d_{\rm IA}(\Phi(b),\pp \Delta_{Q^3})
    = \min_{i=0,\ldots,4} \{\ell_i(\Phi(b))\}. 
\end{equation}
Again, we only discuss the case~$b_1,b_2 > 0$, as the others are analogous. Using \eqref{eq:bigPhi}, we find that the functions~$\ell_i \circ \Phi$ are given by 
\begin{equation*}
    \begin{split}
    \ell_0(\Phi(b)) = \frac{2\pi}{3} &- b_1 - b_2 - 2b_3, \quad
    \ell_1(\Phi(b)) = \frac{2\pi}{3}+b_3,\quad
    \ell_2(\Phi(b)) = \frac{2\pi}{3}+b_1 + b_3,\\
    \ell_3(\Phi(b)) &= \frac{2\pi}{3}+b_1 + b_2 + b_3, \quad
    \ell_4(\Phi(b)) = \frac{2\pi}{3}+b_2 + b_3. 
    \end{split}
\end{equation*}
In that case, we clearly have~$\ell_1 \leqslant \ell_2,\ell_3,\ell_4$, meaning that the minimum in \eqref{eq:PSIcomp} is attained by~$\ell_0$ or by~$\ell_1$. Note that those are indeed two of the expressions appearing in the minimum in \cref{prop:germPchek}. The remaining three cases yield~$b_3$ and the remaining three expressions in the minimum, in each case. 
\noindent
\\
{\bf Step 4: Extension.} The above determines $\psi_{\cp^{\Ch}}$ on an open dense subset of a neighbourhood of the origin. By Proposition~\ref{prop:min_functionals}, this determines $\psi_{\cp^{\Ch}}$ everywhere, and thus $\Delta^{\!\Psi}(\cp^{\Ch})$. This step does not work in the case of displacement energy.
\proofend

\begin{remark}
\label{rk:lagisotopy_11}
In light of this proof, the fact that~$\cp$ and~$\cp^{\Ch}$ are \emph{Lagrangian isotopic} becomes very evident. Indeed, both~$\cp$ and~$\cp^{\Ch}$ can be deformed into Lagrangian tori in a nearby reduced space with~$b_1,b_2 \neq 0$, which only has one orbifold point. Since their curves are isotopic there, so are the respective deformations. 
\end{remark}

\subsection{Lifted Vianna tori}

The quadric $Q^3$ contains infinitely many monotone Lagrangian tori which are not Hamiltonian isotopic (and not symplectomorphic, either) to one another. These can be constructed by lifting Vianna tori from the two-dimensional quadric $Q^2 = S^2 \times S^2$ to $Q^3$ by 
\begin{enumerate}
    \item viewing $Q^2 = S^2 \times S^2$ as a suitable symplectic quotient of $Q^3$ and lifting Vianna tori from that quotient, or, equivalently, by
    \item Biran's circle bundle construction.
\end{enumerate}
This is why we call them \emph{lifted Vianna tori}. Both the construction using Biran's approach and the fact that there are infinitely many tori in $Q^3$ are known to experts, though this result does not seem to have appeared in the literature. For example it can be seen as a variation of the case of $\CP^3$ which is treated in \cite{ChaHirWan24}, or as a special case of forthcoming work by Diogo--Tonkonog--Vianna--Wu \cite{DTVW26}, where similar ideas will be treated in much greater generality. We will use the first approach and work with symplectic reduction.

As we will see below, the Oakley--Usher torus $\cp \subset Q^3$ is the lift of the Clifford torus $S^1_{\rm eq} \times S^1_{\rm eq} \subset S^2 \times S^2$. What is of interest to us here is that the Chekanov Oakley--Usher Lagrangian is not of this form.

\begin{theorem}
    \label{thm:exoticCOU}
    The Chekanov Oakley--Usher Lagrangian $\cp^{\Ch} \subset Q^3$ is not symplectomorphic to any of the lifted Vianna tori. 
\end{theorem}

\begin{remark}
    Note that the case of $\CP^3$ is quite different, since the corresponding tori $\cp_{1,1}, \cp_{1,1}^{\rm Ch} \subset \CP^3$ are Hamiltonian isotopic. However it can be shown that both these tori \emph{are not} symplectomorphic to the ones treated in \cite{ChaHirWan24}. This follows from the computation of Newton polytopes of disk potentials which were carried out in \cite[Section 8]{OakUsh16} and in \cite[Theorem 1.1]{ChaHirWan24}.
\end{remark}

Let us now discuss the Vianna tori $T_k^{\rm Via} \subset S^2 \times S^2$. Let $\Delta_0 = [-\frac{2\pi}{3},\frac{2\pi}{3}] \times [-\frac{2\pi}{3},\frac{2\pi}{3}] \subset \R^2$ be the Delzant polytope of the standard toric structure on $S^2 \times S^2$ equipped with the monotone symplectic form which evaluates to $\frac{4\pi}{3}$ on both standard generators of $H_2(S^2 \times S^2)$. It will become clear below why we choose this normalization. 

Let $\{\Delta_k\}_{k \in \N}$ denote the set of all possible almost toric base diagrams obtained from ATF-mutations of $\Delta_0$. We take $T^{\rm Via}_k$ to be the fibre over $0 \in \Delta_k$, after nodal sliding the nodes close enough to the vertices such that $0$ does not lie on any branch cuts.

\begin{remark} 
Some of the $\Delta_k$ are triangles, and some are quadrilaterals. Here we consider them all, as opposed to \cite{Via17}, where only the triangular base diagrams are considered. The triangular ones have the advantage that they can be enumerated by the solutions of some Markov-type equation, as in the case of $\CP^2$, and thus the set they form has a rather simple description. We again refer to \cite{Via17} for details. The quadrilateral base diagrams on the other hand are ill-understood and, to our knowledge, it is not even known whether all of the corresponding tori can be distinguished. This is not important for us here; we will prove that none of their lifts is the image of $\cp^{\Ch}$ under an ambient symplectomorphism of $Q^3$.
\end{remark}

We can forget all of the almost toric decorations (nodes, branch cuts) of $\Delta_k$ to obtain an ordinary polytope $\Delta_k \subset \R^2$. 

\begin{proposition}
The $\Psi$-polytope of each Vianna torus coincides with its ATF base diagram, 
	\begin{equation}
    	\label{eq:delta_psi_downstairs}
    	\Delta^{\!\Psi}(T^{\rm Via}_k) = \Delta_k. 
	\end{equation}
\end{proposition}
\proof
The proof is based on the same approach as the proof of \cite[Theorem 3.2]{Bre23b}, to which we refer for details. Let $\mu_{S^2\times S^2} \colon S^2 \times S^2 \rightarrow \square$ be the standard toric fibration of $S^2 \times S^2$ over the square $\square = [-a,a] \times [-a,a]$, whose size $a > 0$ depends on the normalization of the symplectic form. For every $y \in \square$, denote by $T_y \subset S^2 \times S^2$ the corresponding toric fibre. Then $\Psi(T_y) = d_{\rm IA}(y,\pp \square)$ by \cite[Proposition C]{SheTonVia24}. Since this agrees with the displacement energy, we can use the same arguments as in the proof of \cite[Theorem 3.2]{Bre23b}, to find that $\Psi(\pi_k^{-1}(y)) = d_{\rm IA}(y,\pp \Delta_k)$, for $y$ in an open dense subset of a neighbourhood of the origin. Here $\pi_k \colon S^2 \times S^2 \rightarrow \Delta_k$ is the almost toric projection map. Thus, $\psi_{T^{\rm Via}_k}$ agrees with the integral affine distance to the boundary of $\Delta_k$ on an open dense subset of a neighborhood of the origin and applying Proposition \ref{prop:min_functionals} again, we find that it does so everywhere. This proves the claim.
\proofend

We note that for the dual of the Newton polytope of the disk counting potential (instead of the $\Psi$-polytope), this result was proved by Vianna \cite[Theorem 5.1]{Via17}. It also follows from subsequent work by Pascaleff--Tonkonog~\cite{PasTon20}. 

Let us now discuss the relevant lifts of these tori using symplectic reduction. To that end, we perform symplectic reduction on $Q^3$ with respect to the third component $\mu_3 \colon Q^3 \rightarrow \R$ of the Gelfand--Cetlin moment map $\mu$. More precisely we consider the symplectic quotient $M = \mu_3^{-1}\left(\frac{4\pi}{3}\right)/S^1$. Since $\mu$ is toric away from $\Sigma$ (see Proposition \ref{prop:GCfibration}), it is easy to check that the $S^1$-action on this level set is free and thus $M$ is an honest symplectic manifold. To identify $M$, we note that it carries a residual toric $T^2$-action which has the polytope 
\begin{equation*}
    \Delta_{Q^3} \cap \left\{x_3 = \textstyle\frac{4\pi}{3}\right\} = \left[0,\textstyle\frac{4\pi}{3}\right] \times \left[0,\textstyle\frac{4\pi}{3}\right] \times \left\{\textstyle\frac{4\pi}{3}\right\}
\end{equation*}
as its image. This polytope is a standard square and thus the reduced space is symplectomorphic to $Q^2 = S^2 \times S^2$ with an appropriately scaled monotone symplectic form. This is an instance of what is called \emph{toric reduction} in \cite[\S 2.2]{Bre23}, to which we refer for details.

Denote by $p \colon \mu_3^{-1}\left(\frac{4\pi}{3}\right) \rightarrow Q^2$ the reduction map. For any Lagrangian $L \subset Q^2$, we denote its lift via $p$ by $\widehat{L} = p^{-1}(L) \subset Q^3$. This lift is automatically a Lagrangian submanifold. 

\begin{proposition}
    The lift $\widehat{T}_{\rm Cl}$ of the Clifford torus $T_{\rm Cl} = S^1_{\rm eq} \times S^1_{\rm eq} \subset M$ coincides with the Oakley--Usher Lagrangian $\cp \subset Q^3$.
\end{proposition}
\proof
    The Clifford torus $T_{\rm Cl} = S^1_{\rm eq} \times S^1_{\rm eq}$ is the fibre over the central point $(\textstyle\frac{2\pi}{3},\textstyle\frac{2\pi}{3}) \in [0,\frac{4\pi}{3}] \times [0,\frac{4\pi}{3}]$ in the moment polytope of $M$. Under the affine-linear isomorphism of this square to $\Delta_{Q^3} \cap \left\{x_3 = \textstyle\frac{4\pi}{3}\right\}$, the point $(\textstyle\frac{2\pi}{3},\textstyle\frac{2\pi}{3})$ maps to the point $(\textstyle\frac{2\pi}{3},\textstyle\frac{2\pi}{3},\textstyle\frac{4\pi}{3}) \in \Delta_{Q^3}$, the fibre of which is precisely $\cp$, which proves the claim. See also \cite[Proposition 2.8]{Bre23} for details. 
\proofend

The $\Psi$-polytope of $\cp$ is (a translate of) $\Delta_{Q^3}$ itself by Corollary \ref{cor:germP}. Since the $\Psi$-polytope of $T_{\rm Cl}$ is the square, this suggests to view $\Delta_{Q^3}$ as the (truncated) cone over the square as follows 
\begin{equation*}
    \Delta_{Q^3}
    = \operatorname{Cone}(S) 
    = \{ tx \sth t \in [0,1], x \in S\},
\end{equation*}
for $S = [0,2\pi] \times [0,2\pi] \times \{2\pi\}$. In other words, the $\Psi$-polytope of the lifted Lagrangian is obtained by taking the cone over the $\Psi$-polytope of the Lagrangian which is lifted.

\begin{lemma}
\label{lem:psi_polytope_lifted}
Let $S_k=(\frac32\Delta_k)\times\{2\pi\}$.
Up to translation, the $\Psi$-polytope of $\widehat T_k^{\rm Via}$ is the truncated cone over $S_k$. More precisely, 
\begin{equation}
\label{eq:psi_polytope_lifted_Vianna}
\Delta^\Psi(\widehat T_k^{\rm Via})+(0,\textstyle{\frac{4\pi}{3}})
=\operatorname{Cone}(S_k).
\end{equation}
In particular, the integral affine distance
from the apex to the opposite facet is $2\pi$.
\end{lemma}

\proof A similar statement was proved in \cite{Bre23b} for the displacement energy of the lift to $\C^3$ of Vianna tori in $\CP^2$. Since the proof is a straightforward adaptation of the proof of \cite[Proposition 1.6]{Bre23b}, we only give an outline here and refer the reader to \cite{Bre23b} for details.

An explicit versal deformation of $\widehat{T}^{\rm Via}_k \subset Q^3$ can be constructed from a versal deformation of the Vianna torus $T^{\rm Via}_k$ downstairs. On an open dense subset, the members of this deformation downstairs are Hamiltonian isotopic to ordinary toric fibres (see \cite[\S 3.3]{Bre23b}). The lifts of these toric fibres are fibres of the Gelfand--Cetlin moment map $\mu \colon Q^3 \rightarrow \Delta_{Q^3}$, of which we know the $\Psi$-invariant. Tracking which member of the versal deformation is Hamiltonian isotopic to which $\mu$-fibre yields the result on an open dense subset of a neighbourhood of the origin. We again use \cref{prop:min_functionals} to deduce that the result holds everywhere. 
\proofend

We can now prove the main theorem of this subsection. \smallskip

\proofof{\cref{thm:exoticCOU}}
It suffices to prove that $\Delta^{\!\Psi}(\cp^{\rm Ch})$ is not integral affine equivalent to any of the $\Delta^{\!\Psi}(\widehat{T}^{\rm Via}_k)$. By Lemma~\ref{lem:psi_polytope_lifted}, we know that the latter $\Psi$-polytopes are cones over the $\Delta_k$, meaning cones over triangles or over quadrilaterals. 
\begin{enumerate}
	\item If $\Delta_k$ is a triangle, then its cone cannot be integral equivalent to $\Delta^{\!\Psi}(\cp^{\rm Ch})$, since the latter contains a four-valent vertex;
	\item If $\Delta_k$ is a quadrilateral, then the integral affine distance of its four-valent vertex to the opposite facet is $2\pi$. The integral affine distance of the four-valent vertex in $\Delta^{\!\Psi}(\cp^{\rm Ch})$ to its opposite facet is $\pi$. Since this is an integral affine invariant, the claim follows.
\end{enumerate}
\proofend

\section{Enumerative properties of the OU-torus in $Q^3$}
\label{sec:enumerative}

In this section, we use a toric degeneration of the quadric $Q^3$ to study the disk-counting potential of the Oakley--Usher torus $\cp \subset Q^3$. Additionally, we apply this to prove a split-generation result for the monotone Fukaya category of $Q^3$ and compare it to Abouzaid--Diogo~\cite{AboDio23}.

\subsection{Superpotentials and non-displaceability}
In \cite{OakUsh16}, the question whether the monotone OU-Lagrangian~$\pkm(r_*)\subset Q^n$ is displaceable for~$k\geqslant 1$ is raised. In this subsection, we give a negative answer for the monotone toric case~$\mathcal{P}= \cp_{1,1}(\frac{4\pi}{3}) \subset Q^3$. To show this, we use a method to reconstruct the superpotential of $\mathcal{P}$ from the versal deformation of its $\Psi$-invariant, and use it to show that this Lagrangian supports critical local systems, in the sense of Biran-Cornea \cite{BirCor12}.
\\
\\
Throughout this section, we work with the universal Novikov field over $\mathbb C$ with real exponents, defined as follows:
\begin{equation}
\label{eq:NovikovField}
  \Lambda
:=
\left\{
\sum_{j=0}^{\infty} a_j T^{\lambda_j}
\ \middle|\
a_j\in \mathbb C,\ \lambda_j\in\mathbb R,\ 
\lambda_j\to +\infty
\right\}.  
\end{equation}
We use multiplicative notation for Novikov weights: a holomorphic disk class $\beta$ contributes the factor $T^{\omega(\beta)}$. All Floer categories and vector-space identifications below use $\Z/2\Z$ gradings.

Note that~$\mathcal{P}$ is a spin manifold: any~$n$-torus is parallelizable (as all Lie groups are), and therefore its second Stiefel--Whitney class~$w_2(\mathcal{P}) \in H^2(\mathcal{P};\Z_2)$ vanishes. Denote by~$\Spin(\mathcal{P})$ the set of its spin structures. Let~$\partial: \pi_2(Q^3,\mathcal{P})\to \pi_1(\mathcal{P})$ be the connecting homomorphism from the homotopy long exact sequence of the pair~$(Q^3,\mathcal{P})$. Let $S \in \Spin(\cp)$. We define the superpotential of~$(\mathcal{P},S)$ to be the map~$\mathscr{W}_{(\mathcal{P}, S)} : \Hom(\pi_1(\mathcal{P}),\Lambda^{\times}) \to \Lambda$,
\begin{equation}
\label{eq:superpotential}
\mathscr{W}_{(\mathcal{P}, S)}(\chi) := \sum_{\substack{A \in \pi_2(Q^3,\mathcal{P})\\ m(A) = 2}} n_{A,S}\,T^{\langle\omega ,A\rangle}\,\chi(\partial A) , 
\end{equation}
where $n_{A,S}\in\Z$ are signed open disk counts: for generic compatible almost complex structures~$J$, they count the number of~$J$-holomorphic disks representing the class~$A$, and whose boundary lies on~$\mathcal{P}$ and passes through a generic point of~$\mathcal{P}$. The choice of spin structure determines the orientations of the moduli spaces of Maslov-2 disks, and hence the coefficients $n_{A,S}$.

\begin{remark}
Two remarks about the definition of the superpotential are in order. 
First, the superpotential is often defined as a map~$\mathscr{W}_{(\mathcal{P}, S)} : \Hom(H_1(\mathcal{P}),\Lambda^{\times}) \to \Lambda$ in the literature. Since $\Lambda^\times$ is abelian, every character
of $\pi_1(\mathcal P)$ factors through
$H_1(\mathcal P)$. Thus the homotopical and
homological definitions of the character torus
agree.

Second, the map~$\partial: \pi_2(Q^3,\mathcal{P})\to \pi_1(\mathcal{P})$ is surjective but not injective (see the exact sequence in \S 2.1 in~\cite{OU1}). However, given~$\gamma \in \pi_1(\mathcal{P})$, there is at most one~$A \in \partial^{-1}(\gamma) \subset \pi_2(Q^3,\mathcal{P})$ such that~$m(A) = 2$. Indeed, two such classes differ by a sphere class $d[\ell]$, whose Maslov index is $6d$; equality of the two indices forces $d=0$. Since $\mathcal P$ is monotone, every Maslov--2
disk has area $2\pi/3$. By Gromov compactness,
only finitely many relative classes contribute
nontrivially to the disk counts. Hence the
superpotential is a Laurent polynomial.

\end{remark}
Choose the basis $\gamma_1,\gamma_2,\gamma_3$ supplied by the effective toric action of $\mu$, as in Corollary~\ref{cor:germP}. The moduli space of rank-one $\Lambda$-valued local systems on~$\mathcal{P}$ is the character torus
\[
\mathcal M(\mathcal{P})
:=
\Hom(\pi_1(\mathcal{P}),\Lambda^{\times})
\cong
(\Lambda^{\times})^3,
\]
where the last identification is given by the holonomy coordinates~$x_i:=\chi(\gamma_i)$. Therefore, in what follows we think of the superpotential as a Laurent polynomial with coefficients in the Novikov field:
\[
\mathscr{W}_{(\mathcal{P}, S)} \in \Lambda[x_1^{\pm 1}, x_2^{\pm 1},x_3^{\pm 1}].
\]
The cohomology ring of~$\mathcal{P}\cong T^3$ is the exterior algebra~$H^*(\mathcal{P};\Z) = \Lambda^*(H^1(\mathcal{P};\Z))$, in particular~$H^*(\mathcal{P};\Z)$ is ring-generated by~$H^1(\mathcal{P};\Z)$ by the usual cup product. Therefore, the Biran--Cornea criterion obtained in \cite[Proposition 3.3.1]{BirCor12} applies: the existence of a spin structure $S$ and a critical point~\(\chi\in (\Lambda^{\times})^3\) of~$\mathscr{W}_{(\mathcal{P}, S)}$, that is, a point satisfying 
\[
\frac{\partial  }{\partial x_i}\mathscr{W}_{(\mathcal{P}, S)}(x_1,x_2,x_3)=0, \qquad i=1,2,3,
\]
implies that the twisted Lagrangian Floer cohomology of the monotone brane~$(\mathcal{P}, S,\chi)$, denoted as~$HF^*((\mathcal{P},\chi),(\mathcal{P},\chi))$, does not vanish, in fact
\begin{equation}
    HF^*((\mathcal{P},\chi),(\mathcal{P},\chi)) \cong H^*(\mathcal{P};\Lambda)
\end{equation}
as $\Z/2\Z$-graded vector spaces. In particular,~$(\mathcal{P},S,\chi)$ is a nonzero object of the monotone Fukaya category, and $\mathcal{P}$ is nondisplaceable as a Lagrangian submanifold. Thus, to prove nondisplaceability, it suffices to choose a spin structure and exhibit a critical point of the corresponding superpotential. In the remainder of the proof, we do exactly this.

\begin{remark}
\label{lemma:support_potential}
	For any pair $(\cp, S)$, we can read the form of the superpotential from the versal deformation of its $\Psi$-invariant. By \cite[Theorem B]{SheTonVia24}, the dual of the Newton polytope of~$\mathscr{W}_{(\mathcal{P}, S)}$ agrees (up to scaling) with 
\begin{equation}
	\label{eq:Psi_polytope_2}
    \Delta^{\!\Psi}(\cp) = \{a \in \R^3 \sth a_1,a_2,a_3-a_1,a_3-a_2,-a_3 \geqslant -\textstyle\frac{2\pi}{3} \}.
\end{equation}
By \eqref{eq:Psi_polytope_2}, the Newton polytope of the potential (with respect to the variables $x_i$ introduced above) is the convex hull of the vectors
\begin{equation}
	\label{eq:exponents}
	(1,0,0),\quad (0,1,0),\quad (-1,0,1),\quad (0,-1,1),\quad (0,0,-1).
\end{equation}
This convex hull does not contain any lattice points besides the vertices and the origin $(0,0,0) \in \Z^3$. The origin would give a constant term  corresponding to a class $\beta\in \pi_2(Q^3,\mathcal{P})$ with $\partial\beta=0$. By exactness, such a class lies in the image of $\pi_2(Q^3)$. Write $j_*(A)=\beta$, for $A \in \pi_2(Q^3)$. Its Maslov index is then
\[
m(\beta)=2c_1(Q^3)(A)\in 2N_{Q^3}\mathbb Z=6\mathbb Z.
\]
Therefore, the constant term vanishes, and the superpotential of the Oakley--Usher torus is of the form
\begin{equation*}
		\mathscr{W}_{(\mathcal{P}, S)}(x_1,x_2,x_3)
			=
		c_1x_1+c_2x_2+c_3\frac{x_3}{x_1}+c_4\frac{x_3}{x_2}+c_5\frac{1}{x_3},
	\end{equation*}
    for non-vanishing coefficients $c_j \in \Lambda^\times$. We compute the full potential in \cref{prop:superpotential of P} using methods which are independent of \cref{sec:tori_quadric}.
\end{remark}
\vspace{1 em}

We now introduce the notion and techniques of toric degenerations. They will be used to compute the superpotential of $\mathcal{P}$ exactly, not just its support. 
\begin{definition}
    A toric degeneration of a symplectic manifold $(X,\omega)$ is a flat family $\Pi: \mathcal X \to \C$, whose fibres $X_t := \Pi^{-1}(t)$ satisfy the following:
    \begin{enumerate}
        \item For $t \neq 0$, the fibre $X_t$ is smooth and symplectomorphic to $X$
        \item $X_0$ is a singular toric variety
        \item The fibres $X_t$ are projective subvarieties of the same projective space, i.e. there is a morphism $\varphi: \mathcal{X}\to \CP^N$, such that for every $t\in \C$, $\varphi|_{X_t}: X_t \to \CP^N$ is an embedding.
    \end{enumerate}
\end{definition}

Toric degenerations have the following useful feature:
symplectic parallel transport along a path in the base
induces symplectomorphisms between the regular fibres
of the degeneration. For a real path, this transport can
be described by the gradient-Hamiltonian flow introduced
by Ruan \cite{Rua01}. When the path approaches the
singular fibre $X_0$, the flow yields a symplectomorphism
away from the vanishing locus, see for example
\cite{NisNohUed10,HarKav15}.

We now fix a toric degeneration $\Pi: \mathcal X \to \C$ by setting
\[
\mathcal X =
\left\{
\bigl([z_0:z_1:z_2:z_3:z_4],t\bigr)\in \mathbb{C}P^4\times \mathbb{C}
\:|\:
tz_0^2+z_1^2+z_2^2+z_3^2+z_4^2=0
\right\},
\]
and taking $\Pi$ to be the projection to the $\C$-component. The fibre $X_1 = \Pi^{-1}(1)$ is by definition just the smooth quadric threefold $Q^3$. Using the gradient-Hamiltonian vector field and a path avoiding the origin, we find that all fibres $X_{t \neq 0}$ are symplectomorphic to $Q^3$. The central fibre is
\[
X_0 =
\left\{
[z_0:z_1:z_2:z_3:z_4]\in \mathbb{C}P^4
\:|\:
z_1^2+z_2^2+z_3^2+z^2_4=0
\right\},
\]
which we can think of as the cone in $\CP^4$ over the two-dimensional quadric in $\{z_0 = 0\} \subset \CP^4$. We will use the gradient-Hamiltonian parallel transport for the path $\gamma \colon [0,1] \rightarrow \C$ defined by $\gamma(t) = 1 - t$. By \cite[Theorem A]{HarKav15}, we find that there is a continuous map $\phi \colon Q^3 = X_1 \rightarrow X_0$ which restricts to a symplectomorphism 
\begin{equation}
    \label{eq:Harada--Kaveh-map}
    \phi\vert_{X_1 \setminus V} \colon X_1 \setminus V \rightarrow X_0 \setminus \{p\},
\end{equation}
where $p=[1:0:0:0:0]$ is the singular point of~$X_0$ and $V = \phi^{-1}(p)$ is the \emph{vanishing locus} of the degeneration $\Pi$ with respect to the path $\gamma$. Let us compute the vanishing locus. Recall from (25) in~\cite{OU1} that we have denoted by $\Sigma \subset Q^3$ the Lagrangian sphere obtained as the image of the zero section in $T^*S^3$ under the symplectic cut. It is explicitly given by $\Sigma = \{ [1:iq] \in \CP^4 \sth q \in S^3\} \subset Q^3$. In the Gelfand--Cetlin fibration introduced in \cref{prop:GCfibration}, it appears as the singular fibre~$\Sigma = \mu^{-1}(0)$ over the origin. 

\begin{lemma}
    \label{lem:vanishingl}
    The vanishing locus $V$ of the degeneration $\Pi$ is $\Sigma \subset Q^3$. 
\end{lemma}

\proof
Equip $\CP^4\times\C$, and thus by restriction
$\mathcal X$, with a product K\"ahler form
\[
\omega_{\mathcal X}
:=
\omega_{FS}\oplus\frac{ic}{2}
dt\wedge d\overline t,\qquad c>0.
\]
We use an antisymplectic involution $\sigma$
on $\mathcal X$ which preserves the fibres
over the real locus and whose fixed locus in
the central fibre $X_0$ consists of the singular
point $p\in X_0$.
The fixed locus is preserved by symplectic
parallel transport along the path $[0,1]\subset\C$.
This strategy is standard; see, for example,
\cite[Lemma 1.20]{Eva24}.

In our setup, the antisymplectic involution is
given by
\[
\sigma:\mathcal X\longrightarrow\mathcal X,
\qquad
([z_0:\ldots:z_4],t)
\longmapsto
([\overline z_0:-\overline z_1:\ldots:
-\overline z_4],\overline t).
\]
Its fixed locus in $X_t$, for real $t>0$, is
\[
\Fix(\sigma)\cap X_t
=
\{[1:i\sqrt t\,q]:q\in S^3\}.
\]
These spheres collapse to $p$ as $t\to0$.
Since $\Fix(\sigma)\cap X_1=\Sigma$, we obtain
$\Sigma\subseteq V$.

To prove the reverse inclusion, consider the
$\SO(4)$-action on the last four homogeneous
coordinates, and let $\mu_{\SO(4)}$ denote
its moment map.
The action preserves the degeneration and
the K\"ahler form. Moreover, its infinitesimal
generators are tangent to the fibres of $\Pi$.
Since the parallel-transport vector field is
symplectically orthogonal to the fibres,
$\mu_{\SO(4)}$ is preserved by parallel transport.

The point $p$ is fixed by $\SO(4)$, so
$\mu_{\SO(4)}(p)=0$.
Consequently,
\[
V\subseteq
\mu_{\SO(4)}^{-1}(0)\cap X_1.
\]
We now show that this set is precisely
$\Sigma$. Indeed, writing $z'=(z_1,\ldots,z_4)=x+iy$,
the moment map vanishes exactly when
$x\wedge y=0$.
In the affine chart $z_0=1$, the equation
of $X_1$ gives
\[
x\cdot y=0,\qquad
|y|^2-|x|^2=1.
\]
Since $x\wedge y=0$, the vectors $x$ and $y$
are linearly dependent. Together with
$x\cdot y=0$, this forces $x=0$ and
$|y|=1$.
Thus the zero set of the moment map in
this chart is exactly
\[
\{[1:iy]:y\in S^3\}=\Sigma.
\]
On the hyperplane $z_0=0$, the quadric
equation instead gives
\[
x\cdot y=0,\qquad |x|^2=|y|^2.
\]
As $x$ and $y$ cannot both vanish, these
conditions are incompatible with
$x\wedge y=0$.
Therefore, there are no additional zeros
of the moment map.

We conclude that
\[
\mu_{\SO(4)}^{-1}(0)\cap X_1=\Sigma,
\]
and hence $V\subseteq\Sigma$. This proves $V=\Sigma$.
\proofend

Let us point out that the antisymplectic involution $\sigma$ is very natural. It extends the antisymplectic involution induced on $Q^3$, after the symplectic cut, by the standard antisymplectic involution $(q,p) \mapsto (q,-p)$ on $T^*S^3$. We are now in a position to prove the first part of \cref{thm:quadric_non_displ} by computing the superpotential of the Oakley--Usher torus. 

\begin{proposition}
\label{prop:superpotential of P}
There exists a reference spin structure \(S_0\in \operatorname{Spin}(\mathcal{P})\) such that the five Maslov--\(2\) disk classes contributing to the superpotential all have algebraic count \(+1\). With respect to this spin structure,
\[
\mathscr{W}_{(\mathcal{P},S_0)}(\chi)
=
T^{2\pi/3}
\left(
x_1+x_2+\frac{x_3}{x_1}+\frac{x_3}{x_2}+\frac{1}{x_3}
\right).
\]
\end{proposition}

\proof
To prove this result, we prove that $\Pi \colon \mathcal{X}
\rightarrow \C$ is a toric degeneration and compute its moment
polytope. Since $\cp$ will turn out to be the monotone fibre of
the toric degeneration, we can apply Nishinou--Nohara--Ueda
\cite[Theorem 1.2]{NisNohUed12} to compute its superpotential.

Recall that by \cite[Theorem A]{HarKav15} and
\cref{lem:vanishingl}, we have found the symplectomorphism
\eqref{eq:Harada--Kaveh-map}. Since
$Q^3\setminus\Sigma=X_1\setminus\Sigma$ is toric by
\cref{prop:GCfibration}, we can transport its toric structure
by the map $\phi$ to find that $X_0\setminus\{p\}$ carries
a symplectic toric structure.

To verify that the central fibre $X_0$ is algebraically toric,
we argue directly. Make a unitary linear change of the four
coordinates $z'=(z_1,z_2,z_3,z_4)$ to coordinates
$w_{00},w_{10},w_{01},w_{11}$ in which its equation is
\[
w_{00}w_{11}=w_{10}w_{01}.
\]
Thus $X_0$ is the projective cone over the Segre embedding
$\CP^1\times\CP^1\subset\CP^3$, and in particular is an
algebraic toric variety. The dense algebraic torus acts with
weights
\[
0,\quad (0,0,1),\quad (1,0,1),\quad
(0,1,1),\quad (1,1,1)
\]
on $z_0,w_{00},w_{10},w_{01},w_{11}$, respectively.
Therefore, with the normalization in which a projective line
has symplectic area $2\pi$, its moment polytope is exactly
the Gelfand--Cetlin polytope $\Delta_{Q^3}$.

It remains to identify the transported torus with the
specific Oakley--Usher torus $\cp$. To this end, observe
that the degeneration is $\SO(4)$-equivariant, where
$\SO(4)$ acts on the last four homogeneous coordinates.
Since the symplectic parallel transport is defined using
an $\SO(4)$-invariant Kähler form, it is
$\SO(4)$-equivariant and preserves the corresponding
moment map.

The two chosen rotation components of this moment map
are precisely the last two components of $\mu_{1,1}$.
Moreover, the norm of the angular-momentum bivector
recovers its radial component. Indeed, in the cut
coordinates
\[
[\sqrt{1-|p|^2}:p+iq],
\]
the appropriately normalized norm is $2\pi|p|$.
On $X_0$, writing $z'=x+iy$, the equation of the
central fibre implies $|x|=|y|$ and $x\cdot y=0$.
Consequently, the same normalized function is
\[
2\pi\frac{|z'|^2}{|z_0|^2+|z'|^2},
\]
which is the toric height coordinate. It follows that, under parallel transport, the effective
toric moment map on $X_0$ agrees with $M\mu_{1,1}$,
up to the fixed choice of signs and ordering of the two
rotation components. In particular, the fibre with
rotation components zero and height $4\pi/3$ is the
image of the Oakley--Usher torus $\cp=\mu^{-1}(x_{\cp})$, as in \cref{prop:GCfibration}. Thus $\phi(\cp)$ is the
monotone toric fibre over $x_{\cp}$ in $X_0$.

The singularity of $X_0$ at $p$ is a \emph{conifold},
or \emph{ordinary double point} singularity. This
follows from the toric normal form of $\Delta_{Q^3}$
at the origin, see \cite[Example 4.32]{Eva23}.
In the same example, Evans discusses why this
singularity admits a small resolution. More explicitly,
a global small resolution of $X_0$ is obtained by
recording one ruling of the rank-one $2\times2$
matrix $(w_{ij})$. The singular point $p$ is replaced
by a copy of $\CP^1$, and no exceptional divisor
is introduced.

Furthermore, $X_0$ is Gorenstein Fano. Indeed, since
it is a normal quadric hypersurface in $\CP^4$, the
adjunction formula gives
\[
K_{X_0}
=
\left(K_{\CP^4}+\mathcal O_{\CP^4}(2)\right)|_{X_0}
=
\mathcal O_{X_0}(-3).
\]
In particular, $K_{X_0}$ is Cartier and
$-K_{X_0}=\mathcal O_{X_0}(3)$ is ample.

We have therefore verified the hypotheses needed to
apply \cite[Theorem 1.2]{NisNohUed12}. We choose the
reference spin structure $S_0$ on $\cp$ by transporting
the standard toric spin structure from the corresponding
toric fibre of $X_0$.

Since the central fibre $X_0$ is toric and admits a small
resolution, the Maslov--$2$ potential of a
Gelfand--Cetlin torus fibre in the smooth fibre is
determined by the toric facets of the polytope associated
to $X_0$. More precisely, if $L(y)=\mu^{-1}(y)$ for
$y\in\operatorname{Int}(\Delta_{Q^3})$, then, with respect
to the reference spin structure $S_0$,
\[
\mathscr W_{(L(y),S_0)}(x)
=
\sum_{j=0}^4 T^{\ell_j(y)}x^{v_j},
\]
where $v_j$ are the primitive inward normals to the five
facets of $\Delta_{Q^3}$, and $\ell_j(y)$ are the
corresponding affine distances to the facets. For the
monotone fibre $\mathcal{P}$, all five affine distances
are equal to $2\pi/3$. The corresponding exponent
vectors are
\[
(1,0,0),\quad (0,1,0),\quad (-1,0,1),\quad
(0,-1,1),\quad (0,0,-1),
\]
from which the desired expression of
$\mathscr W_{(\mathcal{P},S_0)}$ follows.
\proofend

For completeness, we discuss the effect that changing spin structure has on the superpotential. The set~$\Spin(\mathcal{P})$ is an~$H^1(\mathcal{P}; \Z_2)$-torsor: for any given~$S_{0} \in \Spin(\mathcal{P})$, any other spin structure is of the form 
\[
S_0 + \eta, \qquad \eta \in H^1(\mathcal{P}; \Z_2) \cong (\Z_2)^3.
\]
Therefore, there are eight inequivalent spin structures on~$\mathcal{P}$.

\begin{corollary}
    For a general spin structure \(S=S_0+\eta\), with $\eta=(\eta_1,\eta_2,\eta_3)\in H^1(\mathcal P;\mathbb Z_2)$, the superpotential is
\[
\mathscr{W}_{(\mathcal P,S)}(\chi)
=
T^{2\pi/3}
\left(
(-1)^{\eta_1}x_1
+
(-1)^{\eta_2}x_2
+
(-1)^{\eta_1+\eta_3}\frac{x_3}{x_1}
+
(-1)^{\eta_2+\eta_3}\frac{x_3}{x_2}
+
(-1)^{\eta_3}\frac{1}{x_3}
\right).
\]
\end{corollary}
\proof
It follows immediately from the fact that, if we choose the spin structure $S=S_0+\eta$, then the contribution of a disk class $\beta$ is multiplied by $(-1)^{\langle \eta,\partial\beta\rangle}$.
\proofend

For the remainder of the paper, we work with the spin structure $S_0$. For our superpotential~$\mathscr{W}_{(\mathcal{P}, S_{0})}$, we consider its Jacobian ring:
\[
\textnormal{Jac}(\mathscr{W}_{(\mathcal{P}, S_{0})}) := \frac{\Lambda[x_1^{\pm 1}, x_2^{\pm 1},x_3^{\pm 1}]}{\langle x_i\partial_{x_i}\mathscr{W}_{(\mathcal{P}, S_{0})}\:|\: i = 1,2,3\rangle }
\]

\begin{proposition}
\label{prop:critical points}
The critical locus of~\(\mathscr{W}_{(\mathcal{P},S_{0})}\) consists of exactly three nondegenerate points. In
particular,
\begin{equation}
\textnormal{Jac}(\mathscr{W}_{(\mathcal{P},S_{0})})\cong \Lambda\oplus \Lambda\oplus \Lambda
\end{equation}
is semi-simple.
\end{proposition}
\proof
The critical-point equations are
\[
\begin{cases}
x_1-\dfrac{x_3}{x_1}=0,\\[6pt]
x_2-\dfrac{x_3}{x_2}=0,\\[6pt]
\dfrac{x_3}{x_1}+\dfrac{x_3}{x_2}
-\dfrac{1}{x_3}=0.
\end{cases}
\]
The first two equations imply $x_1^2=x_2^2=x_3$,
and hence
\[
(x_1-x_2)(x_1+x_2)=0.
\]
On the other hand, the third equation gives
\[
x_1+x_2=\frac{x_1x_2}{x_3^2}.
\]
Since the right-hand side is invertible in the
Jacobian ring, so is $x_1+x_2$. Consequently,
$x_1=x_2$ holds in the Jacobian ring itself.
Substituting $x_2=x_1$ and $x_3=x_1^2$ into
the third equation yields an isomorphism
\[
\operatorname{Jac}(\mathscr W_{(\mathcal P,S_0)})
\cong
\frac{\Lambda[x_1^{\pm1}]}{\langle x_1^3-\frac12\rangle}.
\]
Since $x_1^3-\frac12$ has three distinct roots
in $\mathbb C\subset\Lambda$, the Chinese remainder
theorem gives
\[
\operatorname{Jac}(\mathscr W_{(\mathcal P,S_0)})
\cong\Lambda\oplus\Lambda\oplus\Lambda.
\]
In particular, the critical scheme consists of
three reduced points, which are therefore
nondegenerate. This also proves that the Jacobian
ring is semisimple.
\proofend

Hence there exist three local systems~$\chi_{\lambda}\in \mathcal{M}(\mathcal{P})$, parametrized by 
\begin{equation}
\label{eq:critical local systems}
\chi_{\lambda} \mapsto (x_1,x_2,x_3) = \left(2^{-1/3}\lambda,2^{-1/3}\lambda,2^{-2/3}\lambda^2\right)\in ( \C^{\times})^3, \qquad \lambda^3=1.
\end{equation}
such that~$HF^*((\cp,\chi_\lambda),(\cp,\chi_\lambda)) \neq 0$. In particular,~\(\mathcal{P}\subset Q^3\) is non-displaceable, proving the first claim in \cref{thm:quadric_non_displ}.

\subsection{Homological mirror symmetry and split-generation}
The preceding superpotential computations may appear ad hoc, but they admit a natural explanation via the mirror symmetry picture. In this subsection, we make this precise by relating our construction to the standard Landau--Ginzburg mirror of the quadric. After that, we show that the OU-Lagrangian $\mathcal{P}$ with its three critical local systems and the Lagrangian sphere $\Sigma$ of Proposition \ref{prop:GCfibration} split-generate the monotone Fukaya category of $Q^3$. Throughout, we denote by $\mathcal F(M)$ the monotone Fukaya category of $M$, and by $D^\pi\mathcal F(M)$ its split-closed derived category.

We begin on the quantum cohomology side (or A-side) of the correspondence. By Beauville's computation of the quantum cohomology of complete intersections
\cite{Bea95}, applied to the quadric hypersurface $Q^3\subset \mathbb CP^4$,
we have 
\begin{equation}
QH^*(Q^3) \cong  \frac{\Lambda[h]}{\langle h^4 -4qh \rangle}, \qquad |h| = 2.
\end{equation}
Here $h \in H^2(Q^3; \Z)$ is the hyperplane class, and $q:=T^{\langle \omega,A \rangle}=T^{2\pi}$ is the quantum parameter associated to the generator $A$ of $H_2(Q^3)$.

We now define the \textit{quantum multiplication operator} $(c_1(Q^3)\,\star -) \in \operatorname{End}(QH^*(Q^3))$, given by 
\[
\alpha \mapsto  c_1(Q^3) \star\alpha = 3(h \star \alpha).
\]
From the ring relation, this operator satisfies $c_1(Q^3)^4=108\,q\,c_1(Q^3)$, so its spectrum is 
\[
\operatorname{Spec}(c_1(Q^3)\,\star -) = \{0\}\sqcup \big\{(108 q)^{1/3} \lambda\: |\: \lambda^3 =1\big\}.
\]
It follows that~$QH^*(Q^3)$ is semi-simple and splits into four eigensummands
\[
QH^*(Q^3) \cong \Lambda e_0 \oplus \bigoplus_{ \lambda^3 =1}\Lambda e_\lambda,
\]
where 
\begin{equation}
    e_0 := 1-\frac{h^3}{4q}, \qquad e_\lambda := \frac{h^3+(4q)^{1/3}\lambda h^2+(4q)^{2/3}\lambda^2 h}{12q}.
\end{equation}
It follows from direct calculation that these eigenvectors of the quantum multiplication operator are idempotent and 
\[
e_0 +\sum_{\lambda^3=1} e_\lambda = 1.
\]

Let us recall the categorical meaning of these eigenvalues and eigenvectors. For a monotone symplectic manifold $M$, its monotone Fukaya category
decomposes as a direct sum of summands indexed by the possible values of $m_0$, the curvature term of the \(A_\infty\)-algebra:
\[
\mathcal F(M)=\bigoplus_{w_\lambda}\mathcal F(M)_{w_\lambda}.
\]
Every summand containing a nonzero object is indexed by an eigenvalue of quantum multiplication by $c_1(M)$. The converse is not automatic for a general monotone $M$; however, it is the case for $M= Q^3$, as we will show below.

 Indeed, for an object \((L,\xi)\), the (length-zero) closed-open string map is a unital algebra homomorphism
\[
\mathcal{CO}^0_{(L,\xi)}:
QH^*(M)
\longrightarrow
HF^*((L,\xi),(L,\xi))
\]
satisfying
\[
\mathcal{CO}^0_{(L,\xi)}(c_1(M))
=
m_0(L,\xi)\,1_{(L,\xi)}.
\]
Thus, if $L$ is a non-trivial object of $\mathcal{F}(M)$ (that is, if $HF^*((L,\xi),(L,\xi))\neq0$), the scalar $m_0(L,\xi)$ must lie in the spectrum of $c_1(M)\,\star -$.

We describe the mirror-symmetry picture in terms of the closed-open string map. The local system chart associated to the monotone torus
\((\mathcal P,S_0)\) gives a Landau--Ginzburg chart with superpotential
\(\mathscr W_{(\mathcal P,S_0)}\). 

On the B-side, the appropriate category is the category of matrix factorizations:
\[
\bigoplus_{\lambda^3=1}\operatorname{MF}\bigl((\Lambda^\times)^3,
\mathscr W_{(\mathcal P,S_0)}-w_\lambda\bigr).
\]
Its three isolated critical points give three Jacobian factors. The category of matrix factorizations of the unshifted $\mathscr W$ on its zero fibre is trivial, since all its critical values are nonzero.

Since $QH^*(Q^3)$ is generated as a $\Lambda$-algebra
by $h=\frac{1}{3} c_1(Q^3)$, the length-zero closed--open map
associated with each critical local system satisfies
\[
\mathcal{CO}^0_{(\cp,\chi_\lambda)}(f(h))
=
f\left(\frac{w_\lambda}{3}\right)\,1_{(\cp,\chi_\lambda)}.
\]
As $HF^*((\cp,\chi_\lambda),(\cp,\chi_\lambda))\neq0$, its unit identifies
the scalar subalgebra
$\Lambda 1_{(\cp,\chi_\lambda)}$
with $\Lambda$. Thus, the images of all elements of
$QH^*(Q^3)$ under these closed--open maps lie in
canonically identified scalar subalgebras.

Composing the three resulting scalar characters with
the inverse of the critical-point evaluation isomorphism
\[
\operatorname{Jac}(W_\cp)
\xrightarrow{\;\cong\;}
\prod_{\lambda^3=1}\Lambda
\]
defines a unital algebra homomorphism
\[
\mathcal{CO}_\cp^0:QH^*(Q^3)\longrightarrow
\operatorname{Jac}(W_\cp).
\]
We refer to $\mathcal{CO}_\cp^0$ as the combined scalar
closed--open homomorphism. It satisfies
\[
\mathcal{CO}_\cp^0(c_1(Q^3))=[W_\cp],
\]
where $[W_\cp]$ is the class of the superpotential
in its Jacobian ring.
\begin{proposition}
    The following holds:
    \begin{itemize}
        \item For all $\lambda,\zeta\in\C$ with~$\lambda^3=\zeta^3=1$, 
        \begin{equation}
            \mathcal{CO}^0_{(\mathcal P,\chi_{\lambda})}(e_{\zeta})=\delta_{\lambda,\zeta}\,1_{(\mathcal P,\chi_{\lambda})},
        \end{equation}
        In particular, $\mathcal{CO}^0_{\mathcal P}$ is surjective.
        \item $\ker(\mathcal{CO}^0_{\mathcal P})=\Lambda\,e_0$.
    \end{itemize}
    Consequently, $\mathcal{CO}^0_{\mathcal P}$ descends to an isomorphism of rings
    \begin{equation}
        QH^*(Q^3)/\Lambda\,e_0 \cong \operatorname{Jac}(\mathscr W_{(\mathcal P,S_0)}).
    \end{equation}
\end{proposition}
\proof
Recall that $\mathcal{CO}^0_{(\mathcal P,\chi_{\lambda})}(h) = (4q)^{1/3}\lambda\,1_{(\mathcal P,\chi_\lambda)}$, therefore 
\[
\mathcal{CO}^0_{(\mathcal P,\chi_{\lambda})}(e_\zeta) = \frac{1}{3}(1+\lambda \zeta^2 +\lambda^2\zeta) = \begin{cases}
    1, \qquad \lambda = \zeta,\\
    0 \:\:\qquad\lambda \neq \zeta.
\end{cases}
\]
The second bullet follows from
\[
\mathcal{CO}^0_{\mathcal P}(e_0) = 1-\frac{\mathcal{CO}^0_{\mathcal P}(h)^3}{4q} = 0.
\]
 As $\mathcal{CO}^0_{\mathcal P}$ is a map of rings, the result follows from the first isomorphism theorem for rings. 
\proofend

We see that the image of the closed-open string map detects exactly the nonzero part of quantum cohomology, while its kernel detects the remaining zero eigensummand. Thus $\mathscr W_{(\mathcal P,S_0)}$ should not be regarded as the full global Landau--Ginzburg mirror potential of $Q^3$, but rather as the potential on the Laurent
chart associated to $\mathcal P$. See \cite{Smi25} for related discussions. With our construction, we can give a simple geometric representative of the missing summand. Namely, the Gelfand--Cetlin fibre over the vertex $(0,0,0)\in \Delta_{Q^3}$ is the Lagrangian sphere
$\Sigma=\mu^{-1}(0)\cong S^3$. The next Proposition shows that the
rank-one summand not seen by the Laurent chart of $\mathcal{P}$ is naturally accounted for by the singular Gelfand--Cetlin sphere fibre, and that these pieces split-generate the monotone Fukaya category of $Q^3$.

\begin{proposition}
\label{prop:generation}
The Lagrangian sphere $\Sigma$ split-generates the summand $D^\pi\mathcal F(Q^3)_0$, and each $(\mathcal{P},\chi_\lambda)$ split-generates the summand $D^\pi\mathcal F(Q^3)_{w_\lambda}$.
Consequently,
\[
D^\pi\mathcal F(Q^3)
=
\langle \Sigma\rangle_{\mathrm{split}}
\oplus
\bigoplus_{\lambda^3=1}
\langle (\mathcal{P},\chi_\lambda)\rangle_{\mathrm{split}}.
\]
\end{proposition}
\proof
We first show that $\Sigma$ defines a nonzero object in the zero eigensummand. 
Since $\pi_1(\Sigma)=\pi_2(\Sigma)=0$, the long exact sequence of the pair gives $\pi_2(Q^3,\Sigma)\cong \pi_2(Q^3)$.
Thus $\Sigma$ has minimal Maslov number $N_{\Sigma}=2N_{Q^3}=6$. In particular, there are no Maslov-2 disks, therefore $m_0(\Sigma)=0$. Since $6>\dim\Sigma+1=4$, Oh's monotone spectral sequence collapses \cite{Oh96}. Thus, as graded vector spaces,
\[
HF^*(\Sigma)\cong H^*(\Sigma; \Lambda)\neq0.
\]
Therefore $\Sigma$ is a nonzero object in the zero eigensummand $\mathcal{F}(Q^3)_0$.

For split-generation, recall that each piece of the semi-simple decomposition of $QH^*(Q^3)$ has rank 1. As $\Sigma$ is a non-trivial object of $\mathcal{F}(Q^3)_0$, and each $(\mathcal{P},\chi_{\lambda})$ is a nontrivial object of $\mathcal{F}(Q^3)_{w_\lambda}$, the result follows immediately from \cite[Corollary 1.12]{She16}.
\proofend

\begin{remark}
This result is the quadric analogue of the Abouzaid--Diogo generation theorem for the compact-monotone Fukaya category of $T^*S^3$ \cite[Theorem 1.4]{AboDio23}. In the next subsection, we recall their results and discuss how the compactification changes the categorical picture.
\end{remark}

\begin{remark}
    This remark was pointed out to us by David Keren Yaar. Since $Q^3$ contains two disjoint monotone spin Lagrangians with non-vanishing Floer cohomology and minimal Maslov number at least two, a result of Alizadeh--Atallah--Cant--Shang \cite[Corollary 21]{AACS25}, applied to $\Sigma$ and $(\cp,\chi_\lambda)$, implies that the spectral diameter is infinite, and hence so is the Hofer diameter of $\Ham(Q^3)$.
 This was already known by other methods \cite[Theorem A]{Kaw24}.
\end{remark}

\begin{corollary}
\label{cor:OU_HMS}
For each $\lambda^3=1$, let
$w_\lambda:=\mathscr W_{(\mathcal P,S_0)}(\chi_\lambda)$.
Then there is an equivalence of split-closed
$\Z/2\Z$-graded triangulated categories
\[
D^\pi\mathcal F(Q^3)_{w_\lambda}
\simeq
D^\pi\operatorname{MF}
\left(
\widehat{\mathcal O}_{(\Lambda^\times)^3,\chi_\lambda},
\widehat{\mathscr W}_{(\mathcal P,S_0)}-w_\lambda
\right),
\]
under which $(\mathcal P,\chi_\lambda)$ corresponds
to the stabilized residue field at $\chi_\lambda$. Moreover, the zero-eigensummand satisfies
\[
D^\pi\mathcal F(Q^3)_0
\simeq
D^\pi(\operatorname{mod}-\mathrm{Cl}_1).
\]
\end{corollary}

\proof
The first assertion follows from
\cite[Theorem C]{Smi25} and
Propositions~\ref{prop:critical points}
and~\ref{prop:generation}.
Indeed, the critical points of
$\mathscr W_{(\mathcal P,S_0)}$ are isolated,
and each corresponding torus brane
split-generates its eigensummand.

The second assertion follows from the description
of the zero-eigensummand of the Fukaya category
of an odd-dimensional quadric in
\cite[Corollary 4.9]{IvanSmith2012}.
\proofend

The preceding corollary gives a geometric realization of
homological mirror symmetry for $Q^3$ in terms of the
Oakley--Usher torus and the Lagrangian sphere $\Sigma$.
Indeed, the three nonzero eigensummands are identified
with the matrix-factorization categories at the three
nondegenerate critical points of the superpotential
$\mathscr W_{(\mathcal P,S_0)}$. The canonical Landau--Ginzburg mirror of $Q^3$, defined in \cite{PechRietschWilliams2016} has an
additional nondegenerate critical point of critical value
zero, which is not contained in the Laurent chart
associated with $\mathcal P$. Its local matrix-factorization
category is graded Morita equivalent to the category
of modules over $\mathrm{Cl}_1$, and hence corresponds
to the zero-eigensummand generated by $\Sigma$. Therefore, the four eigensummands of the Fukaya category
are accounted for by explicit Lagrangian submanifolds.

\subsection{Monotone Fukaya categories under compactification}

At the beginning of this subsection, we recall the setup and split-generation result of \cite[Theorem~1.4]{AboDio23} for the compact-monotone Fukaya category of
$T^*S^n$ (denoted $\mathcal{F}(T^*S^n)$), and its special case $n=3$. After that, we compare with our result for $D^\pi\mathcal{F}(Q^3)$.  We slightly adapt the notation and terminology for consistency with our previous exposition. We denote its split-closed derived category by $D^\pi\mathcal F(T^*S^n)$. Throughout, we assume $n\geq 3$ odd. The even-dimensional case requires different graded modules and bounding-cochain data; see~\cite[Sections~3 and~6]{AboDio23}.

We use the same Novikov field $\Lambda$ as defined in \eqref{eq:NovikovField}. It carries the non-Archimedean valuation $\nu: \Lambda \to \R\cup\{+\infty\}$, defined by
\[
\nu\left(\sum_{j=0}^{\infty}a_jT^{\lambda_j}\right)
=
\min\{\lambda_j\mid a_j\neq 0\},
\qquad
\nu(0)=+\infty.
\]
Following \cite{AboDio23}, we call an element
\emph{unitary} if it is invertible and has zero valuation, and denote the group of unitary Novikov units by
\[
U_\Lambda^*
:=
\{U\in\Lambda^\times\mid \nu(U)=0\}.
\]
Thus a rank-one unitary local system on a Lagrangian $L$ is equivalently
a holonomy homomorphism
\[
\chi:H_1(L)\to U_\Lambda^*.
\]
Let  $F\subset T^*S^n$ be a cotangent fiber, and consider its wrapped $A_\infty$-algebra $\mathcal{A}:= CW^*(F,F;\Lambda)$. The $A_\infty$-Yoneda functor associates to every compact monotone brane
$L$ in $\mathcal{F}$ the right $A_\infty$-module
\[
   L \longrightarrow  CW^*(F,L;\Lambda),
\]
over $\mathcal{A}$. The cohomology algebra of $\mathcal{A}$ is 
\[
A:= H^*(\mathcal A)=HW^*(F,F;\Lambda)\cong H_{-\ast}(\Omega S^n,\Lambda) \cong \Lambda[u], \qquad |u| = 1-n.
\]
Abouzaid--Diogo prove that $A=\Lambda[u]$ is intrinsically formal \cite[Proposition~5.2]{AboDio23}, and moreover that every graded right $A$-module is intrinsically formal
\cite[Proposition~5.4]{AboDio23}. Consequently, cohomology determines the objects of the relevant $A_\infty$-module category up to quasi-isomorphism. We use the cohomological interpretation of \cite[Corollary 5.5]{AboDio23}, with morphisms given by graded Ext groups, throughout the following module discussion.

In particular, the Yoneda image of a
compact monotone brane $L$ may be studied through the ordinary $A$-module
\[
    Y(L):=HW^*(F,L;\Lambda).
\]
Let $\operatorname{mod}_{\mathrm{pr}}(A)$ denote the $\Z/2\Z$-graded cohomological category of finite-dimensional graded $A$-modules, with morphisms $\operatorname{Ext}_A^*(M,N)$ and Yoneda composition. Transport its triangulated structure from the derived $A_\infty$-module category. With this convention the cohomological Yoneda functor is fully faithful, and the generation theorem identifies split-closed triangulated categories:
\[
D^\pi\mathcal F(T^*S^n)\simeq\operatorname{mod}_{\mathrm{pr}}(A).
\]
Thus, the geometric split-generation problem may be studied entirely on
the module side.

Since $\Lambda[u]$ is a PID and $\Lambda$ is algebraically closed, the simple proper graded $\Lambda[u]$-modules, up to parity shift, are precisely 
\[ 
S_a:=\Lambda[u]/\langle u-a \rangle,\qquad a\in\Lambda. 
\]
Every proper $A$-module has finite length and is built from these simple
modules by extensions. Equivalently, the modules $S_a$ generate $\operatorname{mod}_{\mathrm{pr}}(A)$
\cite[Corollary~6.5]{AboDio23}. Therefore, in order to split-generate
$D^\pi\mathcal{F}(T^*S^3)$, it suffices to
realize all the modules $S_a$ geometrically, possibly as direct summands of Yoneda images.

The zero-section contributes the uncountable collection of objects $(0_{S^n},\alpha[\mathrm{pt}])$, where $\alpha[\mathrm{pt}]$ is the family of bounding cochains used in \cite{AboDio23}, with $\nu(\alpha)\geq0$ (including $\alpha=0$). These objects have modules
\begin{equation}
    Y(0_{S^n},\alpha[\mathrm{pt}]) \cong S_\alpha,
\end{equation}
thus producing all the nonnegative-valuation modules. Here and below, we choose the orientation of $[\mathrm{pt}]$
and the sign of $u$ so that the displayed identification holds.

 Second, the negative valuation part is realized by the Oakley--Usher family $\cp_{0,n-1}(r) \subset T^*S^n$, indexed by the area-parameter $r \in (0,+\infty)$. Since $H_1(\mathcal{P}_{0,n-1}(r)) \cong H_1(S^1\times S^{n-1})\cong \mathbb Z$, a rank-one local system is determined by one holonomy $U\in U_\Lambda^*$. Let $\alpha:=T^{-r}\,U^{-1} \in \Lambda$. On the level of $\Lambda[u]$-modules, the brane $(\cp_{0,n-1}(r), U)$ produces the rank-two module
\[
\Lambda[u]/(u^2-\alpha) \cong S_{\sqrt{\alpha}}\oplus S_{-\sqrt{\alpha}}.
\]
Since $U$ is unitary, we have $\nu(\alpha) = -r <0$. Since $n$ is odd, the two summands have the same parity,
and the decomposition is compatible with the chosen
$\Z/2\Z$-grading; see
\cite[Lemma 4.13 and Corollary 4.14]{AboDio23}.

 Thus the zero-section sector covers $\nu(a)\geq 0$, while the $\mathcal{P}_{0,n-1}(r)$-sector covers $\nu(a)<0$. Together they produce all simple modules $S_a$.

For $n=3$, the family $\mathcal{P}_{0,2}(r)$ can be replaced by the torus family $\mathcal{P}_{1,1}(r) \cong T^3$. In a suitable basis of $H_1(T^3) \cong \Z^3$ (described in \cite[Lemma 2.13]{AboDio23}), their superpotential is given by
\begin{equation}
\label{eq:open potential}
\mathscr{W}^{\text{open}}_{\cp_{1,1}(r)}(z_1,z_2,z_3)=T^{r/2}z_1(1+z_2)(1+z_3).
\end{equation}
This superpotential admits critical local systems, given by $\chi_U = (U, -1,-1)$ for $U\in U_\Lambda^*$.
The objects $(\cp_{1,1}(r), \chi_U)$ have modules 
\[
Y(\cp_{1,1}(r),\chi_U)
\cong
S_{\sigma\beta}[\epsilon_1]
\oplus
S_{\sigma\beta}[\epsilon_2],
\]
where $\sigma\in\{1,-1\}$ and
$\epsilon_1,\epsilon_2\in\Z/2\Z$.
Since $\nu(\beta)=-r/2<0$, varying $r>0$
and $U\in U_\Lambda^*$ realizes all the
negative-valuation parameters, independently
of the sign and parity shifts.

We have $\nu(\beta) = -\frac{r}{2}<0$. As $r>0$ and $U\in U_\Lambda^*$ vary, the resulting simple summands, allowing their parity shifts and the fixed sign convention, exhaust all parameters of negative valuation. In this sense, the $S^1\times S^2$-family and the $T^3$-family describe the same negative-valuation part of the module-theoretic picture.

Our compactification result has a different form. In $Q^3$, the monotone Fukaya category decomposes according to the finite spectrum of quantum multiplication by $c_1(Q^3)$:
\[
\operatorname{Spec}(c_1(Q^3)\star-)
=
\{0\}
\sqcup
\{(108q)^{1/3}\lambda\mid \lambda^3=1\}.
\]
The singular Gelfand--Cetlin fibre $\Sigma\simeq S^3$ has
$m_0(\Sigma)=0$ and generates the zero eigensummand. The monotone
Oakley--Usher torus $\mathcal P$, equipped with its three critical local systems
\[
(x_1,x_2,x_3)
=
\bigl(
2^{-1/3}\lambda,\,
2^{-1/3}\lambda,\,
2^{-2/3}\lambda^2
\bigr),
\qquad \lambda^3=1,
\]
generates the three nonzero eigensummands. We compare the torus family before and after compactification. In $Q^3$, only the torus $r = \frac{4\pi}{3}$ is monotone, thus relevant for the monotone Fukaya category $\mathcal{F}(Q^3)$. Indeed, the four disk classes visible in $T^*S^3$ have area $\frac{r}{2}$, while the compactifying disk has area $2\pi-r$, and monotonicity after compactification forces $r=\frac{4\pi}{3}$.

This allows us to study how the superpotential is affected by compactification. For the tori $\cp_{1,1}(r)\subset T^*S^3$, the superpotential is given by \eqref{eq:open potential}. The change of basis 
\[
B=
\begin{pmatrix}
1&0&0\\
-1&1&0\\
-1&-1&1
\end{pmatrix}
\in \GL(3;\mathbb Z)
\]
induces the following monomial change of coordinates on the holonomy variables:
\[(z_1,z_2,z_3) = \left(  x_1,\frac{x_2}{x_1},\frac{x_3}{x_1x_2}\right).
\]
In these variables, the superpotential of $\mathcal{P}$ becomes
\begin{equation*}
    \begin{split}
\mathscr{W}_{\mathcal{P}}(z_1,z_2,z_3) &
=
T^{2\pi/3}
\left(
z_1(1+z_2)(1+z_3)
+
\frac{1}{z_1^2z_2z_3}
\right)\\
& = \mathscr{W}^{\text{open}}_{\cp_{1,1}(4\pi/3)}(z_1,z_2,z_3)+T^{2\pi/3}\frac{1}{z_1^2z_2z_3}.
    \end{split}
\end{equation*}
We see that compactification adds exactly a term to the open superpotential.

 As a result, the critical equations change. The family of unitary critical local systems $\chi_U$ is replaced by three critical local systems, which in these coordinates are 
\[
\widetilde{\chi_\lambda}=(2^{-1/3}\lambda, 1,1), \qquad \lambda^3=1.
\]
Since these coordinates are in $\C^\times \subset \Lambda^\times$, they have valuation zero and thus $\widetilde{\chi_\lambda}$ are still unitary. Moreover, compactification changes the role of the $(S^1\times S^2)$-generators. The monotone element of this family in $Q^3$ is $\mathcal{P}_{0,2}(\frac{4\pi}{3}) \subset Q^3$, and recall from \cite[Proposition~3.16]{OU1} that it is displaceable. Therefore, $\mathcal{P}_{0,2}(\frac{4\pi}{3}) \simeq 0$ in $\mathcal{F}(Q^3)$. Every admissible brane supported on this displaceable compactified Lagrangian has zero Floer cohomology. Thus it contributes no nonzero object after compactification, while suitable branes on the open family contribute to the cotangent-bundle generating collection.

\appendix
\section{Solid mutations (joint with Soham Chanda)}

The goal of this appendix  is to give an alternate description of the Chekanov Oakley--Usher Lagrangian in the quadric $Q^3$ by presenting it as a double \emph{solid torus mutation} of the OU Lagrangian, in the spirit of \cite{ChaHirWan24}. This allows us to compute its disk-potential using wall-crossing results. This requires making a judicious choice of basis for the first homology. To this end, and to simplify our disk-potential computations, we would use the disk-potential computation of \cite{AboDio23}. We again denote by $\cp$ and $\cp^{\rm Ch}$ the Oakley--Usher and the Chekanov Oakley--Usher tori in $Q^3$, respectively.

\subsection{Solid-Torus Mutations and Disk-Potential}

We begin with recalling the notion of solid-torus surgery and how such surgeries affect disk-potentials from \cite{ChaHirWan24}. Geometrically, solid-torus surgery is essentially a circle-fibered version of the disk surgery.

\begin{definition}[Solid Mutation Configuration] A pair $(L,\cT)$ of Lagrangians in a symplectic manifold $M$ is a \emph{solid mutation configuration (SMC)}  if 
    \begin{itemize}
        \item $L$ is a Lagrangian torus;
        \item $\mathcal{T}$ is a Lagrangian solid torus, i.e. $\cT$ is diffeomorphic to $D \times \T^{n-2}$ where $D$ denotes the closed two-dimensional disk ;
        \item $L$ and $\cT$ intersect cleanly along the boundary of $\ct$ 
        \item\label{rigid-torus} 
        the pair $(L,\ct\cap L)$ 
        is diffeomorphic to the standard pair $(\T^n, \T^{n-1})$.
    \end{itemize}
\end{definition}

\begin{example}\label{ex:doubl_torus}
  
From Equation \eqref{eq:OU-param-AD} we see that when $\alpha \in \{ 0, \frac{1}{4} \},$ the image of $I$ restricted to $(0,r) \times \{ \alpha \} \times S^1 \times S^1$ has a solid-torus closure. We denote these solid-tori as $\cT_{0} $ and $ \cT_{\frac{1}{4}}$. We thus have two disjoint solid tori $\cT_0,\cT_{\frac{1}{4}}$ such that $(\cp, \cT_*)$ form a solid mutation configuration. In terms of Figure~2 in~\cite{OU1}, these solid tori project to line segments, at angles $0,\frac{\pi}{4}$ with the horizontal axis, joining the center with the red circle. In coordinates, we have the following presentations, \[
\begin{split}
 \mathcal T_0&=\{q=(x,0),\ p=(0,v):x\in S^1,\ |v|\leq2/3\},\\
 \mathcal T_{1/4}&=\{q=(0,y),\ p=(-v,0):y\in S^1,\ |v|\leq2/3\}.
\end{split}
\]
\end{example}

 An application of a generalized Weinstein neighborhood theorem \cite[Lemma 2.3]{ChaHirWan24}\cite[Theorem A]{KarTukZim26} states that any smooth identification of $L\cup \cT$ as a stratified space can be upgraded to a symplectomorphism from a neighborhood of $L \cup \cT$. This allows one to identify a neighborhood of a solid mutation configuration with neighborhood of a \textit{standard local model} that we describe now. Let $X := \C^2\setminus\{x_1x_2 = 1\}$ and endow it with the Lefschetz fibration $\pi : X \to \C \setminus \{1\}$ defined by $\pi(x) = x_1x_2$. Given a simple loop $\gamma$ in $\C\setminus\{0,1\}$ define the Lagrangian 
\begin{equation}
    T_\gamma:= \left\{ (x,y) \in \C^2  \bigg| |x| = |y|, \; \pi(x,y) \in \im(\gamma)\right\}.
\end{equation}
The standard solid-mutation configuration in $\C^{n}$ is $(T_{S^1} \times \T_{Cl}^{n-2}, T_{[0,1]}\times \T_{Cl}^{n-2})$ where $\T^{n-2}_{Cl}$ is the unit-radius Clifford torus in $\C^{n-2}$.

\begin{definition}
    Let $(L,\cT)$ be a solid mutation configuration and identify a neighborhood $U$ of $L \cap \cT$ to $U_{std}$, a neighborhood of the standard solid mutation configuration,  via a symplectomorphism $\varphi:U \to U_{std}\subset \C^n $ such that $\varphi(L)=T_{S^1} \times \T^{n-2}$ and $\varphi(\cT)=T_{[0,1]}\times \T^{n-2}$. A \textit{solid mutation} for the data $(L,\cT,\varphi)$ replaces $L$ with $L_\mu:=\varphi\inv(T_{\gamma} \times \T^{n-2})$ where $\gamma$ is a curve in $\C \setminus[0,1].$ If $L$ is monotone, one can choose a curve $\gamma$ to ensure $L_\mu$ is monotone. A solid mutation is also called a solid 
torus surgery.
\end{definition}

\noindent When one performs a solid mutation of a monotone torus, the disk-potential change can be explicitly computed. We state the main wall-crossing result needed for our computation without proof.

\begin{proposition} \label{prop:cross}\cite{PasTon20},\cite[Theorem 3.5]{ChaHirWan24}
    Select a basis $(\hat x,\hat y,\hat w_1,\dots \hat w_{n-2})$ of the first homology of $T_{S^1}\times \T^{n-2}$ using the natural splitting $H_1(T_{S^1})\oplus H_1(\T^{n-2})$ and $\hat{x},\hat{y}$ are chosen such that $\hat x$ is the fiber circle in the Lefschetz fibration and $\hat{y}$ projects to a counter-clockwise oriented circle in the base. Then we have the following wall-crossing result for change of disk potentials under solid torus surgery,

    $$W_{L}(x,y,w_1,\dots,w_{n-2}) = W_{L_\mu}(x,y(1+x),w_1,\dots , w_{n-2}).$$
    Here we use $x,y,\dots, w_{n-2}$ to denote the monodromy around $\hat x, \hat y,\dots \hat w_{n-2}$.
\end{proposition}

\subsection{Comparison with the Abouzaid--Diogo model}\label{ss: comparison}
We first relate the coordinates used above to the affine-quadric model
used by Abouzaid--Diogo in \cite[Section 2.2]{AboDio23}. Recall that the
symplectic cut identifies the unit co-disk bundle with the affine chart
$\{Z_0\neq 0\}\subset Q^3$ by
\[
    (q,p)
    \longmapsto
    \left[
        \sqrt{1-|p|^2}:p+iq
    \right].
\]
Thus, writing $w_j=\frac{Z_j}{Z_0}$, we obtain affine coordinates
\[
    w=\frac{p+iq}{\sqrt{1-|p|^2}}.
\]
Since $ Q^3
    =
    \left\{
        [Z_0:\cdots:Z_4]\in\CP^4
        \ \middle|\
        Z_0^2+\cdots+Z_4^2=0
    \right\}$, they satisfy $\sum_{i=1}^4 w_i^2 = -1$. It is convenient to pass to the usual model of affine quadric
\[
    X_3
    =
    \left\{
        (\xi_1,\ldots,\xi_4)\in\C^4
        \ \middle|\
        \xi_1^2+\cdots+\xi_4^2=1
    \right\}
\]
by setting $\xi_j:=-iw_j.$
Abouzaid--Diogo use the equivalent affine model
\[
    X_{\mathrm{AD}}
    =
    \left\{
    (z,u_1,v_1,u_2,v_2)\in\C^5
    \ \middle|\
    u_1v_1=z+1,\quad
    u_2v_2=z-1
    \right\}.
\]
The identification between $X_3$ and $X_{\mathrm{AD}}$ is given by
\begin{equation}
\label{eq:quadric-AD-change}
\begin{aligned}
    u_1 &= \sqrt{2}\,(\xi_1+i\xi_2),
    &
    v_1 &= \sqrt{2}\,(\xi_1-i\xi_2),\\
    u_2 &= \sqrt{2}\,(\xi_3+i\xi_4),
    &
    v_2 &= -\sqrt{2}\,(\xi_3-i\xi_4),\\
    z
    & =
    2(\xi_1^2+\xi_2^2)-1
    =
    1-2(\xi_3^2+\xi_4^2).
\end{aligned}
\end{equation}
This is an isomorphism of affine quadrics.

The projection
\[
    z\colon X_{\mathrm{AD}}\longrightarrow\C
\]
is the Morse--Bott Lefschetz-type fibration considered by
Abouzaid--Diogo. Its critical values are $\pm1$. For a regular value
$z\neq\pm1$, they consider the two vanishing circles
\[
    V_{1,z}
    =
    \left\{
        u_1v_1=z+1,\ |u_1|=|v_1|
    \right\},
    \qquad 
    V_{2,z}
    =
    \left\{
        u_2v_2=z-1,\ |u_2|=|v_2|
    \right\}.
\]

\noindent For any curve $C\subset \C$, the Lagrangian  $T_C \subset X_{AD}$ is defined by

$$T_C:= \bigcup_{z\in C} V_{1,z}\times V_{2,z}.$$

\begin{proposition}\label{prop:comp-coordinates}
    We have the following identifications 
    \[
        \cp = T_{C_{\rm OU}}, \quad
        \cT_{0} = T_{[1,\frac{13}{5}]}, \quad
        \cT_{\frac{1}{4}} = T_{[-\frac{13}{5},-1]},
    \]
    where $C_{\rm OU} \subset \C$ is an ellipse. 
\end{proposition}

\proof 
This is a straightforward computation using the above.
\proofend

\begin{remark}[Homology conventions]\label{rem:homconv}
A presentation of disk-potential needs a choice of a basis of first homology. We recall the choices necessary for the calculation of disk-potential as done in \cite{AboDio23,CPU16}.
In particular, we will construct a basis $\langle \hat z_1, \hat z_2, \hat z_3\rangle$ of $H_1(\cp)$ such that
setting the monodromy dual to the basis as $z_*$, the disk potential for the Oakley-Usher torus is 
$$W_\cp (z_1,z_2,z_3) = z_1(z_2+1)(z_3+1) + \frac{1}{z_1^2z_2z_3}.$$

We fix notation for the homology classes represented by the circles $V_{*,z}$ as $$\hat z_2= [V_{2,z}], \,\hat z_3 = [V_{1,z}].$$ We further let $\hat z_1$ denote the unique element that projects to the counter-clockwise orientation of $C_{OU}$ and forms a basis $\langle \hat z_1, \hat z_2 ,\hat z_3 \rangle$ of $H_1(\cp).$

Note that  $\langle \hat z_2, \hat z_3 \rangle$ generates the homology of a fiber $T_z \subset T_{C_{OU}}$. Moreover, the computations for Proposition \ref{prop:comp} show that in $\ct_0$ the $V_{2,z}$ circle collapses and $\ct_{\frac{1}{4}}$ the $V_{1,z}$ circle collapses. Thus
we have the following homological relations from the inclusion of the intersection of the solid tori and $\cp$,

$$ \Z\langle \hat z_2,\hat z_3 \rangle = H_1(\ct_0 \cap \cp) \xrightarrow{\iota^{0}_*} H_1(\ct_0) \quad ( \hat z_2,\hat z_3) \mapsto \hat z_3 $$
$$ \Z\langle \hat z_2,\hat z_3 \rangle = H_1(\ct_{\frac{1}{4}} \cap \cp) \xrightarrow{\iota^{1/4}_*} H_1(\ct_0) \quad ( \hat z_2,\hat z_3) \mapsto \hat z_2. $$
\noindent Notice that from the homology conventions needed for wall-crossing result of Proposition \ref{prop:cross}, we have 
$$\text{mutation along } \ct_{0}\quad\leftrightarrow \quad z_1\mapsto \frac{z_1}{z_2+1}.$$ Similarly $$\text{mutation along } \ct_{\frac{1}{4}}\quad\leftrightarrow \quad z_1\mapsto \frac{z_1}{z_3+1}.$$ 
\end{remark}

\begin{proposition}\label{prop:comp}
     The Oakley-Usher torus $\cp$ and the solid tori $\cT_{0},\cT_{\frac{1}{4}
    }$ are equal to the Abouzaid-Diogo Lagrangian tori $T_{C_{OU}}$, $T_{[1,\frac{13}{5}]}$ and $T_{[-\frac{13}{5},-1]}$ where $C_{OU}$ is an ellipse. Moreover, the Chekanov Oakley-Usher torus $\cp^{\Ch}$ is equal to a double solid torus surgery of $T_{C_{OU}}$ along the solid tori $T_{[1,\frac{13}{5}]}$ and $T_{[-\frac{13}{5},-1]}$.
\end{proposition}

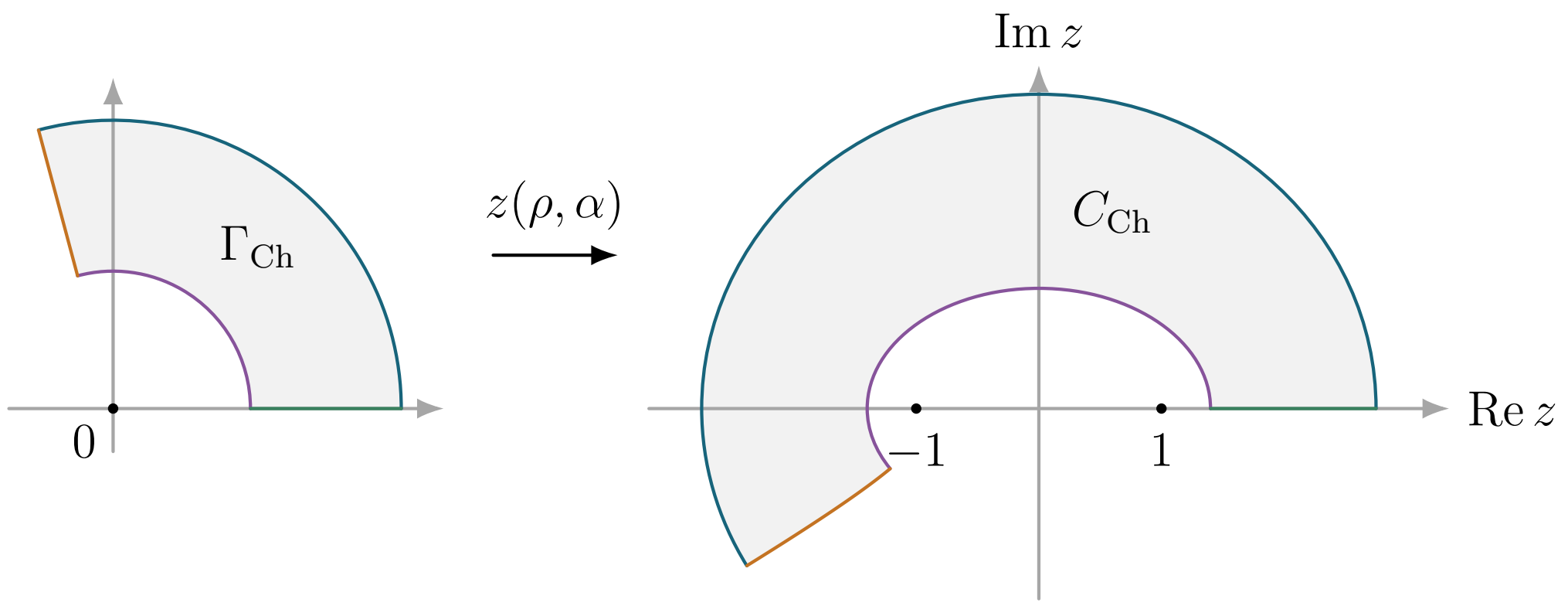
\begin{figure}[htbp]
    \centering
    \begin{tikzpicture}[scale=1.05,>=Latex,
      line cap=round,line join=round]

      \fill[gray!10]
        (2.35,0)
        arc[start angle=0,end angle=105,radius=2.35]
        --(105:1.12)
        arc[start angle=105,end angle=0,radius=1.12]
        --cycle;

      \draw[->,thick,gray!70] (-.85,0)--(2.7,0);
      \draw[->,thick,gray!70] (0,-.35)--(0,2.7);

      \draw[AppOuter,thick]
        (2.35,0)
        arc[start angle=0,end angle=105,radius=2.35];
      \draw[AppInner,thick]
        (1.12,0)
        arc[start angle=0,end angle=105,radius=1.12];
      \draw[AppZero,thick] (1.12,0)--(2.35,0);
      \draw[AppQuarter,thick] (105:1.12)--(105:2.35);

      \fill (0,0) circle (1.2pt);
      \node[below left] at (0,0) {$0$};
      \node at (1.18,1.32) {$\Gamma_{\rm Ch}$};

      \draw[->,thick] (3.1,1.25)--(4.12,1.25);
      \node[above] at (3.6,1.32) {$z(\rho,\alpha)$};

      \begin{scope}[shift={(7.55,0)}]
        \fill[gray!10]
          (2.75,0)
          arc[start angle=0,end angle=210,
              x radius=2.75,y radius=2.56173769]
          -- plot[domain=2.75:1.4,samples=60,variable=\s]
             ({\s*cos(210)},{sqrt(\s*\s-1)*sin(210)})
          arc[start angle=210,end angle=0,
              x radius=1.4,y radius=.97979590]
          --cycle;

        \draw[->,thick,gray!70] (-3.18,0)--(3.35,0)
          node[right,black] {$\operatorname{Re}z$};
        \draw[->,thick,gray!70] (0,-1.55)--(0,2.8)
          node[above,black] {$\operatorname{Im}z$};

        \draw[AppOuter,thick]
          (2.75,0)
          arc[start angle=0,end angle=210,
              x radius=2.75,y radius=2.56173769];
        \draw[AppInner,thick]
          (1.4,0)
          arc[start angle=0,end angle=210,
              x radius=1.4,y radius=.97979590];

        \draw[AppQuarter,thick]
          plot[domain=1.4:2.75,samples=60,variable=\s]
            ({\s*cos(210)},{sqrt(\s*\s-1)*sin(210)});
        \draw[AppZero,thick] (1.4,0)--(2.75,0);

        \fill (-1,0) circle (1.2pt);
        \fill (1,0) circle (1.2pt);
        \node[below] at (-1,-.07) {$-1$};
        \node[below] at (1,-.07) {$1$};
        \node at (.6,1.6) {$C_{\rm Ch}$};
      \end{scope}

    \end{tikzpicture}
    \caption{The curve $\Gamma_{\rm Ch}$ in the Brendel--Piccirilli
      model, drawn in polar coordinates $w=\rho e^{2\pi i\alpha}$
      (left), and its image $C_{\rm Ch}$ in the Abouzaid--Diogo base
      under the coordinate change \eqref{eq:OU-AD-base-alpha} (right).
      Matching colours identify corresponding pieces: the teal and
      purple arcs have constant $\rho$, while the green and orange
      pieces have constant $\alpha$. The critical values $\pm1$ lie
      outside the region bounded by $C_{\rm Ch}$.
      The curves are schematic, with corners understood to be smoothed.}
    \label{fig:cou-in-ad}
\end{figure}

\begin{proof} We start with directly checking that the Oakley--Usher parametrization is
adapted to precisely the circles $V_{\ast,z}$. Let $ x=(x_1,x_2),\; y=(y_1,y_2) \in S^1\subset\R^2$. For $k=m=1$, the Oakley--Usher
parametrization takes the form
\begin{equation}
\label{eq:OU-param-AD}
    q=(\cos(2\pi \alpha)x,\sin(2\pi \alpha)y),
    \qquad
    p=(-\rho\sin(2 \pi \alpha)x,\rho\cos(2\pi \alpha)y).
\end{equation}
In particular, $ |p|=\rho$. Set $X:=x_1+ix_2$, $ Y:=y_1+iy_2\in S^1$, 
and 
\[
    a(\rho,t)
    :=
    \frac{\cos(2\pi\alpha)+i\rho\sin(2\pi\alpha)}
         {\sqrt{1-\rho^2}},
    \qquad
    b(\rho,t)
    :=
    \frac{\sin(2\pi\alpha)-i\rho\cos(2\pi\alpha)}
         {\sqrt{1-\rho^2}}.
\]
Then \eqref{eq:quadric-AD-change} becomes
\begin{equation}
\label{eq:OU-AD-coordinates}
\begin{aligned}
    u_1 &= \sqrt{2}\,a(\rho,\alpha)X,
    &
    v_1 &= \sqrt{2}\,a(\rho,\alpha)X^{-1},\\
    u_2 &= \sqrt{2}\,b(\rho,\alpha)Y,
    &
    v_2 &= -\sqrt{2}\,b(\rho,\alpha)Y^{-1}.
\end{aligned}
\end{equation}
In particular, $|u_1|=|v_1|$ and $|u_2|=|v_2|$. Hence, for fixed $(\rho,\alpha)$, varying $(x,y)\in S^1\times S^1$
parametrizes exactly $V_{1,z(\rho,\alpha)}\times V_{2,z(\rho,\alpha)}$ inside the fibre of the Abouzaid--Diogo fibration.

It remains to determine the corresponding point in the base. By direct computation, we obtain

\begin{equation}
\label{eq:OU-AD-base-alpha}
    z(\rho,\alpha)
    =
    \frac{1+\rho^2}{1-\rho^2}\cos(4\pi\alpha)
    +
    i\frac{2\rho}{1-\rho^2}\sin(4\pi\alpha).
\end{equation}
For fixed $\rho$, the Oakley--Usher Lagrangian therefore projects to
the ellipse
\[
    C_\rho
    :=
    \left\{
        \frac{1+\rho^2}{1-\rho^2}\cos(4\pi\alpha)
        +
        i\frac{2\rho}{1-\rho^2}\sin(4\pi\alpha)
        \ \middle|\
        \alpha \in\R/\Z
    \right\}.
\]
For the monotone value $\rho=2/3$, this is $ C_{\mathrm{OU}}
    :=
    \left\{
        \frac{13}{5}\cos(4\pi\alpha)
        +
        i\frac{12}{5}\sin(4\pi\alpha)
    \right\}$.
Thus the monotone Oakley--Usher torus is precisely the torus obtained
by taking $V_{1,z}\times V_{2,z}$ over the ellipse $C_{\mathrm{OU}}$.

We can now identify the two solid mutation configurations.
The first solid torus is given by
\[
    \cT_0 = S_1
    =
    \left\{
        \alpha=0,\quad
        0\leq\rho\leq\frac23
    \right\}.
\]
It is immediately checked that $z(S_1)
    =
    \left[1,\frac{13}{5}\right]$. Similarly, the second solid torus is
\[
    \cT_{\frac{1}{4}} = S_2
    =
    \left\{
        \alpha=\frac{1}{4},\quad
        0\leq\rho\leq\frac23
    \right\}.
\]
Its image under the projection is $z(S_2)
    =
    \left[-\frac{13}{5},-1\right]$.

We have therefore identified the two solid mutation configurations
with the two natural thimbles of the Abouzaid--Diogo fibration cut off
by the monotone Oakley--Usher ellipse:
\begin{equation}\label{eq:solid} 
    S_1
    =
    \bigcup_{z\in[1,13/5]}
    V_{1,z}\times V_{2,z},
    \qquad
    S_2
    =
    \bigcup_{z\in[-13/5,-1]}
    V_{1,z}\times V_{2,z},
\end{equation}
where at $z=1$ the second circle collapses and at $z=-1$ the first
circle collapses.

On the other hand, the base change equation \eqref{eq:OU-AD-base-alpha} is a diffeomorphism from $(0,1)\times\bigl(\mathbb R/(\tfrac12\mathbb Z)\bigr)$ to $\C\setminus [-1,1]$. A direct check shows that the curve defining the Chekanov Oakley-Usher torus is mapped to a curve $C_{\Ch}$ in the Abouzaid-Diogo base that does not contain the line segment $[-1,1]$ in its interior; see Figure \ref{fig:cou-in-ad}.
\end{proof}

The last claim of the Proposition follows by noticing that the curve $C_{\Ch}$ can be Hamiltonian isotoped, with respect to a correct choice of symplectic form such that Hamiltonian isotopy of base induces Hamiltonian isotopy of the Lagrangian torus above it, to a curve that is realized by a double solid tori surgery of $T_{C_{OU}}$ along the solid tori $T_{[1,\frac{13}{5}]}$ and $T_{[-\frac{13}{5},-1]}$.

\subsection{Disk-Potential Computation}

\begin{theorem}
    The disk potential of the Chekanov OU torus is $$W_{\cp^{\Ch}}(z_1,z_2,z_3) = z_1 + \frac{(z_2+1)^2(z_3+1)^2}{z_1^2z_2z_3}.$$
\end{theorem}

\proof
    From Example \ref{ex:doubl_torus} we know that $(\cp, \ct_0)$ and $(\cp,\cT_{\frac{1}{4}})$ form a solid mutation configuration.  
    Moreover, from an application of the generalized Weinstein neighborhood theorem \cite{KarTukZim26} and Figure~2 in~\cite{OU1}, one can check that the Chekanov OU torus $\cp^{\Ch}$ is the iterated solid torus mutation of $\cp$ along the solid tori $\ct_{0}$ and $\ct_{\frac{1}{4}}$. 
    
$$\cp \xlongrightarrow{\ct_{0}-\text{mutation}} \cp_{\rm Kim}\xlongrightarrow{\ct_{1/4}-\text{mutation}} \cp^{\Ch} $$
    The disk-potential of $\cp$ can be written, using the homology conventions set in Remark \ref{rem:homconv}, as 

    $$W_{\cp}=z_1(z_2+1)(z_3+1) + \frac{1}{z_1^2z_2z_3}.$$

    \noindent

    Using the mutation formula from Proposition \ref{prop:cross} and homology identifications in Remark \ref{rem:homconv}, we have
    $$W_{\cp^{\Ch}}(z_1,z_2,z_3) = W_{\cp}\left(\frac{z_1}{(z_2+1)(z_3+1)},z_2,z_3 \right) = z_1 + \frac{(z_2+1)^2(z_3+1)^2}{z_1^2z_2z_3}. $$
\proofend

\begin{remark}[Kim's Torus]
    We note that the torus obtained from one solid mutation of $\cp$ can be shown to be same as the lift of the Chekanov torus of the two-dimensional quadric sitting as a divisor in $Q^3$. This can be proved by creating affine coordinates near the divisor $Q^2$ where the Lagrangian solid torus in $Q^3$ projects to a Lagrangian disk. In such a local model, one can directly verify that a solid torus surgery corresponds to a lifting of a disk surgery on the base.
 The wall-crossing computation for a single solid torus mutation matches the disk-potential computation of \cite{Kim23}.
\end{remark}

\bibliographystyle{abbrv}
\bibliography{OU2}

\end{document}